\documentclass[12pt,a4paper]{article}
\usepackage[utf8]{inputenc}
\usepackage[margin=1.5cm]{geometry}

\usepackage{bm} 
\usepackage{enumitem} 
\usepackage{array}  
\usepackage{graphicx} 
\usepackage{caption}
\usepackage{amssymb} 
\usepackage{amsmath} 
\usepackage[hidelinks]{hyperref} 
\hypersetup{colorlinks=True,
            citecolor=blue} 
\usepackage{xcolor} 
\usepackage{multirow} 
\usepackage{stmaryrd} 
\usepackage{mathtools} 
\usepackage{authblk} 

\usepackage[subpreambles=true]{standalone} 
\usepackage{tikz}     

\usepackage[style = apa, uniquename=false]{biblatex}
\DeclareLanguageMapping{english}{english-apa}
\usepackage{amsthm}
\newtheorem{theorem}{Theorem}
\newtheorem{lemma}[theorem]{Lemma}
\newtheorem{corollary}[theorem]{Corollary}
\theoremstyle{definition}
\newtheorem{definition}{Definition}

\newenvironment{pf}[1][Proof]{\begin{proof}[#1]}{\end{proof}}

\newcommand{\pfrac}[2]{\frac{\partial{#1}}{\partial{#2}}}
\newcommand{\dx}[1]{\hspace{0.75mm}\mathrm{d}{#1}}
\newcommand{\der}[2]{\frac{\dx{#1}}{\dx{#2}}}

\newcommand{\intOmega}[1]{\int_{\Omega} #1 ~ \text{d}\bm{x}}
\newcommand{\intGamma}[1]{\int_{\Gamma} #1 ~ \text{d}S}
\newcommand{\bracfrac}[2]{\left( \frac{#1}{#2} \right)}
\newcommand{\bracked}[1]{\left( #1 \right)} 
\newcommand{\jump}[1]{\left\llbracket #1 \right\rrbracket_+}

\let\vec\mathbf

\title{A conservative, discontinuous Galerkin, tracer transport scheme using compatible finite elements}
\author[1,2]{Timothy C. Andrews \thanks{timand@umich.edu}}
\author[3]{Thomas M. Bendall}
\affil[1]{Department of Mathematics and Statistics, University of Exeter, Exeter, United Kingdom}
\affil[2]{Department of Climate and Space Sciences and Engineering, University of Michigan, Ann Arbor, MI, USA}
\affil[3]{Dynamical and ML Modelling, Met Office, Met Office, Exeter, United Kingdom}

\begin{document}

\maketitle

\begin{abstract}
This paper outlines a conservative transport scheme for scalar tracers within a compatible finite element model for geophysical fluid equations. Instead of using the advective transport equation for a tracer mixing ratio, the transport equation is cast into conservative form that includes the density of a reference field. This ensures mass conservation in the continuous equations, which can be preserved in the discrete equations with a discontinuous Galerkin transport scheme. Our method is designed to work for two placements of the mixing ratio in a Charney-Phillips vertical staggering: either co-located with the reference density or vertically staggered from it. The new scheme is designed to conserve tracer mass and ensure consistency by maintaining a constant mixing ratio. Additionally, a mass-conserving limiter is developed to enforce monotonicity in the co-located configuration, whilst preserving linear correlations. A number of tests are presented, including with terminator toy chemistry to test the inclusion of physics terms, and a moist rising bubble to test the conservation of moisture species in a dynamical core setup.
\end{abstract}

\section{Introduction}
\label{section:intro}
Mass conservation of dry air and other scalar tracers is one of the key properties of an effective dynamical core for numerical weather prediction and longer climate simulations \parencite{Staniforth_hor_grids,thuburn2008some}. Accordingly, many design decisions in operational and research dynamical cores, such as the choice of transport scheme, are made with a focus on mass conservation. Some transport schemes are inherently mass-conserving, such as a finite volume discretisation of a flux-form conservative equation. Such a discretisation ensures local mass conservation --- that mass lost from one grid cell is transferred to adjacent cells --- which implies mass conservation over the entire domain. Methods that are not generally mass-conserving, such as semi-Lagrangian transport, can often be reformulated in a conservative manner, e.g. \textcite{zerroukat2002slice} and \textcite{lauritzen2010conservative}. Other approaches may be globally but not locally mass-conserving, whilst other methods that are not globally conservative typically use mass-fixers to compensate. This work will focus on the discontinuous Galerkin (DG) transport method, which uses a weak form of the governing equations and permits discontinuities between elements. Like the finite volume method, DG discretisations are locally, and thus globally, conservative if the continuous equations are written in a conservative (flux) form. DG methods are a popular choice for numerical fluid models as they are inherently parallelisable and well-suited for capturing discontinuous behaviours, such as shocks \parencite{nair2011emerging}. \par
The transport component of a dynamical core determines the evolution of a field due to the winds. Such transported fields include the potential temperature and dry air density. Additionally, any number of passive or active tracers may also be transported, including aerosols and moisture variables such as water vapour and cloud water fractions. Tracer species are typically represented by a mixing ratio, $m$, which is the ratio of the \textit{tracer density}, $\rho_X$, to a reference density, $\rho$:
\begin{equation}
m := \rho_X / \rho.
\label{eq:mixing_ratio}
\end{equation}

\noindent To ensure mass conservation of a species during transport, we require a scheme for the reference density and mixing ratio such that the globally integrated tracer density remains constant. Note that this work will take the reference density to be the density of dry air, i.e. the density of air excluding tracers, although all the methods presented here generalise to any choice of reference density. \par
There are two types of scalar transport terms that appear in the governing equations used in dynamical cores; these are the \textit{conservative} and \textit{advective} forms. The dry density obeys a conservative transport equation,
\begin{equation}
    \pfrac{\rho}{t} + \bm{\nabla} \cdot (\rho \vec{u}) = 0,
\label{eq:conservative_rho_transport}
\end{equation}

\noindent with $\vec{u}$ the transporting velocity field, whilst mixing ratios are governed by the advective form,
\begin{equation}
    \pfrac{m}{t} + (\vec{u} \cdot \bm{\nabla}) m = 0.
\label{eq:advective_m_X_transport}
\end{equation}

\noindent Whilst a DG discretisation of the conservative transport equation (\ref{eq:conservative_rho_transport}) will be mass-conserving, this is not the case for the advective form (\ref{eq:advective_m_X_transport}), as discrete updates for the mixing ratio using advective fluxes do not ensure conservation of the tracer density. To resolve this, we instead consider the time evolution of the tracer density, $\rho_X = \rho m$, which for smooth fields is
\begin{equation}
    \pfrac{}{t} (\rho m) = \pfrac{\rho}{t} m + \rho \pfrac{m}{t} = - m \bm{\nabla} \cdot (\rho \vec{u}) - \rho (\vec{u} \cdot \bm{\nabla}) m,
\end{equation}

\noindent using (\ref{eq:conservative_rho_transport}) and (\ref{eq:advective_m_X_transport}), which leads to the conservative equation
\begin{equation}
    \pfrac{}{t} (\rho m) + \bm{\nabla} \cdot
     (\rho m \vec{u}) = 0.
\label{eq:conservative_m_X_transport}
\end{equation}

\noindent We term (\ref{eq:conservative_m_X_transport}) a \textit{conservative tracer} transport equation to distinguish it from the conservative transport equation for the dry density \eqref{eq:conservative_rho_transport}. A DG discretisation of (\ref{eq:conservative_m_X_transport}) will enforce mass conservation as the discrete updates apply to $\rho m$. Most dynamical cores use the conservative tracer equation in preference to an advective form transport of $m$, with some models even using tracer densities as prognostic variables in place of mixing ratios, e.g. MPAS \parencite{skamarock2012multiscale} and NICAM \parencite{satoh2008nonhydrostatic}. \par
Whilst DG methods for conservative tracer transport have been explored previously, such as in \textcite{giraldo1997lagrange}, \textcite{nair2005discontinuous}, and \textcite{bosler2019conservative}, the novelty of this work is that we allow $\rho$ and $m$ to be defined in different function spaces. The use of different function spaces is a key part of the compatible finite element method, which carefully selects spaces that ensure mimetic properties, specifically that the discretisation obeys curl(grad) = $\vec{0}$ and div(curl) = 0 \parencite{cotter_shipton_mixed_FE,cotter2023_acta_numerica,natale2016compatible}. GungHo, the dynamical core for the next-generation LFRic-Atmosphere model from the UK Met Office \parencite{adams2019lfric,melvin2019cartesian,melvin2024mixed}, uses a compatible finite element spatial discretisation. However, GungHo does not solve the transport terms with a finite element method, but instead opts for the flux-form semi-Lagrangian transport scheme, SWIFT \parencite{bendall2025swift}. \par
GungHo uses a Charney-Phillips vertical staggering \parencite{melvin2018choice}, which offsets the dry density from thermodynamic variables such as the potential temperature. Advantages of the Charney-Phillips grid include more accurate gravity wave dispersion and avoiding the computational mode of the Lorenz grid \parencite{thuburn2005vertical}; a drawback is that energy conservation is harder to achieve \parencite{lee2025helmholtz}. The Charney-Phillips configuration offers a choice of mixing ratio placement: (1) $m$ could be co-located with $\rho$ which directly enables mass conservation or (2) $m$ could be co-located with the thermodynamic variables to better capture latent heat associated with phase changes, and radiation/cloud interactions \parencite{konor2000choice,bendall2023trilemma}.  The second placement means that the mixing ratio and dry density are vertically staggered. Given the arguments for either placement of $m$, we design a transport scheme that works for both configurations. \par
With a vertical staggering of $\rho$ and $m$, the mixing ratio lives in a function space that is continuous in the vertical dimension, so we use the embedded DG scheme of \textcite{cotter2016embedded} to map $m$ to a fully discontinuous transport space. After transport, the fields are mapped back to their original spaces using Galerkin projections. As projecting the mixing ratio in isolation will violate mass conservation, we introduce a conservative Galerkin projection that acts on the tracer density. We note that our scheme is only globally conservative in the vertically staggered case, compared to local conservation for the co-located case; local conservation could be enforced in the vertically staggered case by introducing a shifted mesh \parencite{bendall2023trilemma}. \par
In addition to ensuring mass conservation, the new conservative tracer transport scheme is designed to obey the \textit{consistency} property, which requires that a spatially constant mixing ratio remains constant after transport. This is a desirable property for a transport scheme, as when $m$ is constant, the conservative tracer equation (\ref{eq:conservative_m_X_transport}) reduces to a scaled conservative equation for the dry density (\ref{eq:conservative_rho_transport}). Hence, any variation in the constant mixing ratio over time implies an inconsistent numerical treatment of the conservative equations for $\rho_X$ and $\rho$ \parencite{zhang2008consistency,machenhauer2009finite}. The last key property of our scheme is monotonicity, with the challenge being to retain mass conservation after applying the monotonic limiter. Our approach is to blend the higher-order nonmonotonic field with a lower-order representation that is guaranteed to be monotonic. Given the complexity of ensuring mass-conserving limiting with partially continuous function spaces, we focus on a monotonic limiter for the co-located configuration, with the extension to the vertically staggered case being future work. \par
Our conservative tracer transport scheme will be tested on meshes comprised of quadrilateral elements, although the method generalises to other meshes. The co-location of $\rho$ and $m$ is tested on a spherical domain with a gnomonic cubed-sphere mesh. The vertical staggering of $\rho$ and $m$ uses a vertical slice domain, with elements constructed as the tensor product of one-dimensional spaces along the horizontal and vertical dimensions. A range of numerical tests will be performed using the Gusto code library \parencite{shipton2026gusto}, which is built upon the Firedrake finite element codebase \parencite{FiredrakeUserManual}. Gusto solves geophysical fluid equation sets using compatible finite elements and is designed for rapid prototyping. \par
The next section will provide important information about the new transport scheme, including the weak form DG discretisations and the operators used to map fields between function spaces. We provide proofs of conservation and consistency for the transport and mapping operators. We also describe the monotonic limiter for the co-located configuration. Section \ref{section:colocate_spaces} tests the conservation and consistency properties of the scheme with co-located spaces, as well as monotonicity and the preservation of linear correlations. Section \ref{section:vert_stag_spaces} then tests conservation and consistency with a vertical staggering of $\rho$ and $m$, including in a dynamical core-like simulation. Concluding remarks are given in Section \ref{section:conclusions}. \par

\section{Background}
\label{section:background}

This section details important mathematical background and discussion of the conservative transport scheme, including proofs of conservation and consistency. 

\subsection[Discrete function spaces]{Discrete function spaces}
\label{subsec:function_spaces}
The discrete finite element function spaces used in this work are constructed from discontinuous Lagrange, $\dx{Q}_k$, or continuous Lagrange, $Q_k$, elements. The discontinuous Lagrange space is spanned by $ k$th-order polynomials. The continuous Lagrange space is a subset of $\dx{Q}_k$ with the additional constraint of $\mathcal{C}^0$ continuity between elements, which requires nodes on the facets. Two-dimensional function spaces are formed by the tensor product $\otimes$ of two one-dimensional function spaces (i.e. intervals) \parencite{mcrae2016automated}, defined as
\begin{equation}
    \mathbb{V}_A \otimes \mathbb{V}_B = \text{span} \{ f \cdot g | f\in \mathbb{V}_A , g \in \mathbb{V}_B\}.
\end{equation}

In this work, we use next-to-lowest order ($k=1$) function spaces for the density and mixing ratio. Our transport scheme could also be applied with lowest-order ($k=0$) elements by using the recovery scheme of \textcite{bendall2019recovered}, which recovers a lowest order field into a next-to-lowest order space for transport. \par
\begin{figure}
    \centering
    \includegraphics[width=0.6\textwidth]{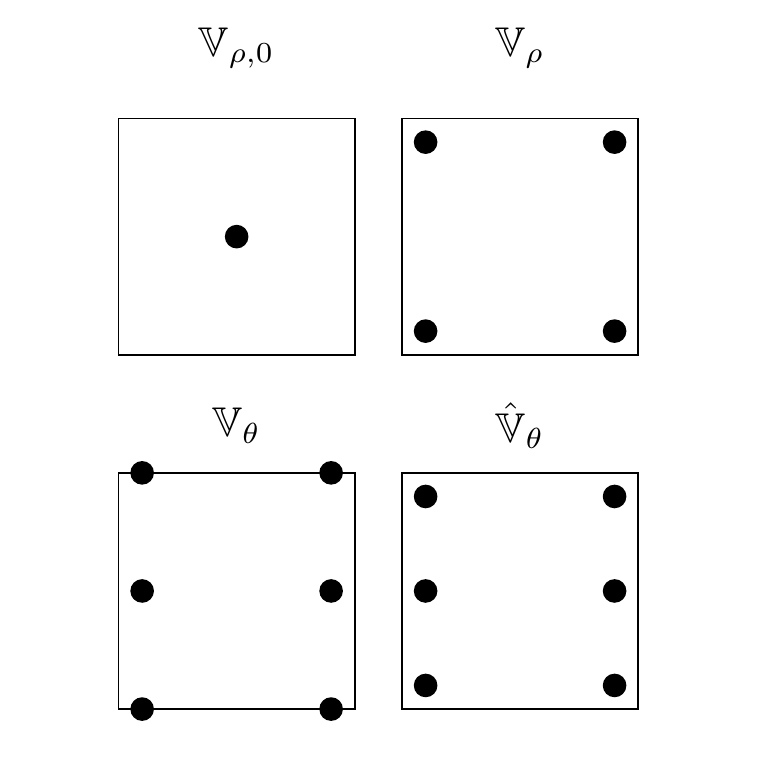}
    \caption{The discrete finite element function spaces used in this work, with the black circles denoting nodes. $\mathbb{V}_{\theta}$ is formed on an extruded mesh that is continuous in the vertical dimension and discontinuous in the horizontal dimension.}
    \label{fig:function_spaces}
\end{figure}

We consider two function spaces for the $k=1$ density and mixing ratio fields (Figure \ref{fig:function_spaces}). The first is the fully discontinuous function space $\mathbb{V}_{\rho}:= \dx{Q}_1 \otimes \dx{Q}_1$. We also consider the thermodynamic function space in a Charney-Phillips vertical staggering, which is continuous in the vertical dimension but discontinuous in the horizontal dimension. The next-to-lowest order thermodynamic function space is $\mathbb{V}_{\theta}:= \dx{Q}_1 \otimes {Q}_2$. We consider two setups: i) the co-located case, where $\rho \in \mathbb{V}_{\rho}, m \in \mathbb{V}_{\rho}$, and ii) the vertically staggered case where $\rho \in \mathbb{V}_{\rho}, m \in \mathbb{V}_{\theta}$. \par
Also shown in Figure \ref{fig:function_spaces} are two other discrete function spaces used in this work. $\mathbb{V}_{\rho,0}$ is the lowest order fully discontinuous space; it is used to compute cell mean fields, such as the cell mean tracer density, and the lowest-order field used in the monotonic limiter. $\hat{\mathbb{V}}_{\theta} := \dx{Q}_1 \otimes \dx{Q}_2$ is a fully discontinuous version of the thermodynamic space $\mathbb{V}_{\theta}$ and is used to transport the prognostic fields in the vertically staggered setup. \par

\subsection{Scheme properties}
\label{subsec:scheme_props}
We now provide a mathematical description of our conservative transport scheme.  Throughout, we require that the density is strictly positive, including after transport and mapping between function spaces,
\begin{equation}
    \rho(\bm{x}) > 0, \quad \forall \bm{x} \in \Omega,
\label{eq:positive_rho}
\end{equation}

\noindent where $\Omega$ denotes the domain and $\bm{x}$ the spatial coordinates. The domain for the co-located case is the surface of the sphere, which is free of boundaries. For the vertically staggered case, the domain is a vertical slice  with a horizontal-vertical ($x$, $z$) Cartesian coordinate system. The boundaries in the horizontal direction are periodic, whilst the boundaries in the vertical dimension enforce impermeability of the flow, $\vec{u} \cdot \widehat{\bm{n}}=0$, for velocity $\vec{u}$ and unit normal at the surface $\widehat{\bm{n}}$. The thermodynamic $\mathbb{V}_{\theta}$ space (Figure \ref{fig:function_spaces}), used for the mixing ratio in the vertically staggered case, is discontinuous in $x$ and continuous in $z$. \par
The conservative tracer transport scheme is designed to ensure the key properties of conservation and consistency. Additionally, in the co-located case, we enforce monotonicity. We now formally introduce these properties.

\begin{definition}[Conservation]
    A conservative scheme ensures that the total tracer mass $M$ is constant in time,
    \begin{equation}
        \der{M}{t} = 0, \quad M = \intOmega{\rho_X},
        \label{eq:dM_dt}
    \end{equation}

    \noindent with $\rho_X = \rho m$ the tracer density.
\end{definition}

 In a discrete model, conservation requires that the total mass is unchanged between timesteps, so
\begin{equation}
        \intOmega{\rho^{n+1} m^{n+1}} = \intOmega{\rho^{n} m^{n}},
\end{equation}
\noindent where superscripts denote the time index. We compute $\rho_X$ as a mean value in each cell by solving the Galerkin projection
\begin{equation}
    \intOmega{\psi \rho_X} = \intOmega{\psi \rho m}, \quad \forall \psi \in \mathbb{V}_{\rho,0},
\end{equation}

\noindent where $\rho$ and $m$ may lie in different function spaces. \par
The numerical tests will either compute the relative change in tracer mass over time,
\begin{equation}
    \Delta M (t) = \frac{| M(t) - M(0)| }{M(0)},
\label{eq:td_error_over_time}
\end{equation}

\noindent with $|\cdot|$ the absolute value and $M(0)$ the initial mass, or a mean mass error over $N$ time samples,
\begin{equation}
    \overline{\Delta M} = \frac{1}{N} \sum_{i=1}^{N} \Delta M (t_i).
\label{eq:mean_td_error}
\end{equation}

\begin{definition}[Consistency]
    A consistent scheme ensures that a spatially constant mixing ratio remains unchanged over time.
\end{definition}

Consistency requires that $m^n=C \rightarrow m^{n+1} = C$, for some spatially constant function $C$.  In numerical tests, this is checked through the relative $L_2$ error
\begin{equation}
    \Delta m(t) = \frac{||m(t) - m(0)||_2}{||m(0)||_2},
\label{eq:m_L2_err}
\end{equation}
\noindent where $|| \cdot ||_2$ denotes the $L_2$ norm, and $m(0)$ is the initial mixing ratio field, which is constant in the consistency tests.

\begin{definition}[Monotonicity]
    A monotonic scheme ensures that no new extrema are produced in the mixing ratio after transport. For DG upwind, this means that the field remains bounded by previous extrema in that element and upwind neighbouring facets.
\end{definition}

\subsection{Formulation for conservative tracer transport}
We now present the mixed finite element equations that describe the transport of the density and mixing ratio. Although $\rho$ and $m$ may live in different function spaces, we ensure that they are represented in a single transport function space of $\mathbb{V}_T$ before performing the transport. For the vertically staggered setup, this requires an appropriate mapping of $\rho, m$ from their original space into the transport space, and then back to their original space after transport. In both the co-located and vertically staggered cases, the transport space $\mathbb{V}_T$ is fully discontinuous: $\mathbb{V}_T = \mathbb{V}_{\rho}$ for the co-located case and $\mathbb{V}_T = \mathbb{\hat{V}}_{\theta}$ for the vertically staggered case (Figure \ref{fig:function_spaces}). The transporting velocity $\vec{u}$ is defined in the Raviart-Thomas function space, which ensures continuity of normal components on each facet. For most of Section \ref{section:background}, we consider the transport of a single mixing ratio, but the scheme generalises to an arbitrary number of mixing ratios. We will later test the transport of multiple mixing ratios in both the co-located and vertically staggered configurations. \par
Consider $\rho, m \in \mathbb{V}_{T}$, with a conservative transport equation \eqref{eq:conservative_rho_transport} for $\rho$ and a conservative tracer transport equation \eqref{eq:conservative_m_X_transport} for $m$. These equations are multiplied by independent scalar test functions $\gamma, \eta \in \mathbb{V}_{T}$ and integrated over the domain $\Omega$ to obtain the weak forms
\begin{subequations}
\begin{align}
 \intOmega{\gamma\pfrac{\rho}{t}} + \intOmega{\gamma\bm{\nabla\cdot} \left(\rho \bm{u}\right)} &= 0, \quad \forall \gamma\in\mathbb{V}_T, \\
 \intOmega{\eta\pfrac{}{t} (\rho m)} + \intOmega{ \eta \bm{\nabla\cdot}\left( \rho m \vec{u} \right)} &= 0, \quad \forall \eta\in\mathbb{V}_T. \label{eq:tracer_conservative_continuous}
\end{align}
\end{subequations}

\noindent Next, we apply integration by parts and the divergence theorem, which leaves surface integral terms defined on the set of interior facets, $\Gamma$. The surface terms are evaluated using discontinuous Galerkin upwinding, leaving
\begin{subequations} \label{eq:mixed_dg_upwind_conserve}
\begin{align}
 \intOmega{\gamma\pfrac{\rho}{t}} + \intGamma{\left(\bm{u} \bm{\cdot}\widehat{\bm{n}}^+\right) \jump{\gamma} \rho^\dagger } - \intOmega{\rho \left(\bm{u}\bm{\cdot\nabla}\right)\gamma} &= 0, \quad \forall \gamma \in\mathbb{V}_T, 
\label{eq:rho_dg_upwind} \\
 \intOmega{\eta \pfrac{(\rho m)}{t}} + \intGamma{\left(\bm{u} \bm{\cdot} \widehat{\bm{n}}^+\right) \jump{\eta} \rho^\dagger m^\dagger }
 - \intOmega{\rho m \left(\bm{u}\bm{\cdot\nabla}\right)\eta} &= 0, \quad \forall \eta\in\mathbb{V}_T, \label{eq:conservative_m_dg_upwind}
\end{align}
\end{subequations}

\noindent with the two sides of a facet arbitrarily labelled as `+' and `-', so $\widehat{\bm{n}}^+$ is an outward facing normal vector on the `+' side of the facet. The jump of a scalar field across a facet is defined as $\jump{\gamma}:= \gamma^+ - \gamma^-$. The dagger $^\dagger$ denotes the upwind value at a facet, e.g.
\begin{equation} \label{def:upwind}
\rho^\dagger := \left\lbrace 
\begin{matrix}
\rho^+ & \mathrm{if} \ \bm{u \cdot}\widehat{\bm{n}}^+ \geq 0, \\
\rho^- & \mathrm{if} \ \bm{u \cdot}\widehat{\bm{n}}^+ < 0.
\end{matrix}\right.
\end{equation} 

\noindent Note that in general, (\ref{eq:mixed_dg_upwind_conserve}) may include inflow and outflow boundary terms \parencite{gibson2019compatible}, but these are not present in the spherical or vertical slice domains we consider. \par
For completeness, the DG upwind formulation of an advective transport equation for the mixing ratio (\ref{eq:advective_m_X_transport}) is
\begin{equation}
    \intOmega{\chi\pfrac{m}{t}} + \intGamma{\left(\bm{u} \bm{\cdot} \widehat{\bm{n}}^+ \right) \jump{\chi}  m^\dagger} - \intOmega{m \bm{\nabla} \cdot (\chi \vec{u})} = 0, \quad \forall \chi\in\mathbb{V}_T,
\label{eq:advective_m_dg_upwind}
\end{equation}

\noindent with $\chi$ a scalar test function. The advective scheme, which solves (\ref{eq:rho_dg_upwind}) and (\ref{eq:advective_m_dg_upwind}), will be compared with the conservative scheme of (\ref{eq:rho_dg_upwind}) and (\ref{eq:conservative_m_dg_upwind}) in the numerical tests of Sections \ref{section:colocate_spaces} and \ref{section:vert_stag_spaces}. \par

\subsection{Time discretisation}
We apply an explicit Runge-Kutta (RK) time discretisation to the DG upwinded transport equations \parencite{cockburn2001runge}. Specifically, we use the third-order strong stability preserving scheme, SSPRK3 \parencite{shu1988efficient}. Strong stability preserving RK methods can be expressed as a sequence of forward Euler steps, allowing for higher-order accuracy whilst retaining the strong stability properties of a single forward Euler step (Section 2.3.3 of \textcite{Durran}). Using a multistage time discretisation like an explicit RK method also simplifies the application of limiters \parencite{cotter2016embedded}. The combination of $k=1$ finite elements and third-order SSPRK timestepping leads to a maximum allowable Courant number of 0.409 \parencite{cockburn2001runge}. \par
The three-stage SSPRK3 scheme applied to a variable $q$ is given by
\begin{subequations}
    \begin{align}
        q^1 &= \mathcal{A} (q^n), \label{eq:ssprk3_stage1}\\
        q^2 &= \frac{3}{4} q^n + \frac{1}{4} \mathcal{A} (q^1), \label{eq:ssprk3_stage2} \\
        \mathcal{T} q^n = q^{n+1} &= \frac{1}{3} q^n + \frac{2}{3} \mathcal{A} (q^2), \label{eq:ssprk3_stage3}
    \end{align}
\label{eq:ssprk3}
\end{subequations}

\noindent where $\mathcal{A}(q)=\widetilde{q}$ represents a forward Euler update from the transport equation, and we denote the updated field with a tilde. Superscripts such as $i$ and $n$ denote the time or stage index. Each SSPRK3 stage can be expressed in the general form
\begin{equation}
    q^{i} = \alpha_i q^n + (1-\alpha_i) \mathcal{A}(q^{i-1}),
\label{eq:ssprk3_gen_stage}
\end{equation}

\noindent where $i \in \{1,2,3\}$ is the stage index, $\alpha_1=0, \alpha_2=3/4, \alpha_3=1/3$, $q^0=q^n$, and $q^3 = q^{n+1}$. We represent one full application of SSPRK3 from time level $n$ to $n+1$ as $\mathcal{T}:\mathbb{V}_T \rightarrow \mathbb{V}_T$. \par
 Our approach is to solve directly for $\rho$ and $m$ at each stage of the conservative transport scheme; solving for $m$ rather than the tracer density of $\rho_X = \rho m$ is necessary for the application of a monotonic limiter on the mixing ratio (Section \ref{subsec:monotonicity}). For the advective scheme, the transport equations for the density and mixing ratios can be combined into a single residual which is solved using a linear solver from the PETSc library \parencite{balay1997efficient,balay2024petsc}. However, for the conservative transport scheme, the combined residual contains nonlinear $\rho m$ terms that would require the use of an iterative nonlinear solver. We instead perform two direct linear solves: we first solve the density equation to find $\rho^{i}$, update the conservative tracer equation with the new density, then solve for the updated mixing ratio $m^i$. Practically, the second linear solve is for the SSPRK3 stage applied to $q=\rho m$, with $\rho^{i} m^{i}$ terms that are linear in the unknown $m^{i}$. \par
We now prove that each forward Euler transport update using the conservative transport system is conservative and consistent. We consider proofs with a combined residual for the density and mixing ratio transport equations and assume that the discrete equations are solved exactly.

\begin{theorem}[Conservative transport]
\label{theorem:conservation}
    The forward Euler update $\mathcal{A}$ of the conservative transport system conserves tracer mass.
\end{theorem}

\begin{pf}
    Combine weak transport equations for the density (\ref{eq:rho_dg_upwind}) and tracer density (\ref{eq:conservative_m_dg_upwind}), and apply a forward Euler discretisation of the time derivative,
    \begin{align}
        &\intOmega{\gamma \widetilde{\rho}} - \intOmega{\gamma \rho} + \Delta t\intGamma{\left(\bm{u} \bm{\cdot}\widehat{\bm{n}}^+\right) \jump{\gamma} \rho^\dagger } - \Delta t \intOmega{\rho \left(\bm{u}\bm{\cdot\nabla}\right)\gamma} + \nonumber \\ &\intOmega{\eta \widetilde{\rho} \widetilde{m}} - \intOmega{\eta \rho m} + \Delta t \intGamma{\left(\bm{u} \bm{\cdot} \widehat{\bm{n}}^+\right) \jump{\eta} \rho^\dagger m^\dagger} - \Delta t \intOmega{\rho m \left(\bm{u}\bm{\cdot\nabla}\right)\eta} = 0, \quad \forall \gamma, \eta \in\mathbb{V}_T,
    \label{eq:theorem_1_residual}
    \end{align}
    \noindent where $\widetilde{\rho}, \widetilde{m}$ are the updated fields after transport. Since our weak form of the transport equation is true for all test functions $\gamma,\eta\in\mathbb{V}_T$, it is true for the choice $\gamma=1, \eta=0$, which gives
    \begin{equation}
        \intOmega{\widetilde{\rho}} = \intOmega{\rho},
    \end{equation}
    \noindent as the jump and gradient terms vanish with this choice of test functions. It is also true that (\ref{eq:theorem_1_residual}) holds for $\gamma=0, \eta=1$, which shows tracer mass conservation of
    \begin{equation}
        \intOmega{\widetilde{\rho}\widetilde{m}} = \intOmega{\rho m}.
    \end{equation} 
\end{pf}

\begin{theorem}[Consistent transport]
\label{theorem:consistency}
    The forward Euler update $\mathcal{A}$ of the conservative transport system is consistent, provided that the same quadrature rule is used for the density and conservative tracer transport equations.
\end{theorem}

\begin{pf}
    Let the mixing ratio be a spatially constant field, $m=C$. Consider the tracer density equation (\ref{eq:conservative_m_dg_upwind}) and subtract $C$ times the density equation (\ref{eq:rho_dg_upwind}), 
    \begin{align}
        &\intOmega{\eta \widetilde{\rho} \widetilde{m}} - \intOmega{\eta \rho C} + \Delta t \intGamma{\left(\bm{u} \bm{\cdot} \widehat{\bm{n}}^+\right) \jump{\eta} \rho^\dagger C^\dagger} - \Delta t \intOmega{\rho C \left(\bm{u}\bm{\cdot\nabla}\right)\eta}  \\  &- \intOmega{C \gamma \widetilde{\rho}} + \intOmega{C \gamma \rho} - \Delta t\intGamma{ C  \left(\bm{u} \bm{\cdot}\widehat{\bm{n}}^+\right) \jump{\gamma} \rho^\dagger } + \Delta t \intOmega{C \rho \left(\bm{u}\bm{\cdot\nabla}\right)\gamma} = 0, \quad \forall \gamma, \eta \in\mathbb{V}_T.
    \end{align}
    We set $\gamma = \eta$, which leads to the cancellation of most terms when the density and tracer density transport equations use the same quadrature rule. This leaves
    \begin{align}
        \intOmega{\eta \widetilde{\rho} \widetilde{m}} -\intOmega{C \eta \widetilde{\rho}} &= 0, \\
        \rightarrow \intOmega{\eta \widetilde{\rho} (\widetilde{m}-C)} &= 0.
    \end{align}
    
    \noindent Choosing $\eta = \widetilde{m}-C$ gives
    \begin{equation}
        \intOmega{\widetilde{\rho} (\widetilde{m}-C)^2} = 0,
    \end{equation}
    which, given that $\widetilde{\rho}>0$ (\ref{eq:positive_rho}), requires that $\widetilde{m}=C$. 
\end{pf}

\begin{corollary}
\label{corollary:ssprk_conservative_consistent}
    Each SSPRK3 stage of the conservative transport scheme is conservative and consistent.
\end{corollary}

\begin{pf}
    From (\ref{eq:ssprk3_gen_stage}), the density and tracer density in each SSPRK3 stage are the convex combination of the time level $n$ field and a forward Euler update. At each stage, the forward Euler update from the previous stage (time level $n$ for stage 1) is conservative and consistent by Theorems \ref{theorem:conservation} and \ref{theorem:consistency}. It thus follows that the solution at each stage is also conservative and consistent. 
\end{pf}

For the co-located case, transport is performed directly in the $\mathbb{V}_{\rho}$ space that $\rho, m$ live in, so the conservative scheme will be conservative and consistent by Corollary \ref{corollary:ssprk_conservative_consistent}. For the vertically staggered case, $\rho$ and $m$ lie in different function spaces and a few more steps are required to ensure conservation and consistency, which are now described. \par

\subsection{Application to vertically staggered spaces}
With the vertical staggering of $\rho \in \mathbb{V}_{\rho}, m \in \mathbb{V}_{\theta}$, the prognostic fields are mapped to the shared $\mathbb{V}_T=\mathbb{\hat{V}}_{\theta}$ space for transport, then projected back.  Our method is based on the embedded DG approach of \textcite{cotter2016embedded}, with the novelty that our injection and projection operators are defined to be mass conserving and consistent. The projection operations are linear in the mapped field and use direct linear solves with an LU factorisation. The general approach for mapping fields between function spaces is the Galerkin projection.

\begin{definition}[Galerkin projection]
    The Galerkin projection $\mathcal{P}: \mathbb{V}^O \rightarrow \mathbb{V}^M$ takes a field in the original function space, $q^O \in \mathbb{V}^O$, and identifies the field in a mapped function space, $q^M \in \mathbb{V}^M$, by solving
\begin{equation}
    \intOmega{\psi q^M} = \intOmega{\psi q^O}, ~\forall \psi \in \mathbb{V}^M.
\label{eq:galerkin_proj}
\end{equation}
\end{definition}

The Galerkin projection (\ref{eq:galerkin_proj}) ensures global conservation of $q$, as seen by choosing $\psi=1$ in (\ref{eq:galerkin_proj}), and hence is sufficient for projecting the density. However, applying $\mathcal{P}$ to a mixing ratio does not ensure that tracer mass is conserved. Accordingly, we introduce a conservative version of the Galerkin projection operator.

\begin{definition}[Conservative Galerkin projection]
    The conservative Galerkin projection of $\mathcal{P}_{\rho m}: \mathbb{V}_{m}^O, \mathbb{V}_{\rho}^O, \mathbb{V}_{\rho}^M \rightarrow \mathbb{V}_{m}^M$ computes $m^M\in\mathbb{V}_m^M$, given fields $m^O\in\mathbb{V}_{m}^O, ~\rho^O \in \mathbb{V}_{\rho}^O, ~\rho^M\in\mathbb{V}_{\rho}^M$, by solving
\begin{equation}
    \intOmega{ \psi \rho^M m^M} = \intOmega{ \psi \rho^O m^O }, \quad \forall \psi \in \mathbb{V}_m^M.
\label{eq:conservative_proj_not_consistent}
\end{equation}
\end{definition}

The conservative Galerkin projection is applied after the projection of the density, $\mathcal{P}:\mathbb{V}_{\rho}^O \rightarrow \mathbb{V}_{\rho}^M$, and has a mass matrix that is weighted by $\rho^M$. As we assume strictly positive density fields (\ref{eq:positive_rho}), the mass matrix is nonsingular. We also assume the density is not small, so that the matrix is not ill-conditioned.  \par

\begin{lemma}
\label{lemma:conservative_projection}
    The conservative Galerkin projection operator is conservative. Furthermore, it is consistent when $\rho^M$ and $m^M$ lie in the same function space and $\rho^M$ is the Galerkin projection of $\rho^O$.
\end{lemma}

\begin{pf}
    For conservation, set $\psi=1$ in the conservative Galerkin equation (\ref{eq:conservative_proj_not_consistent}). \par
    For consistency, set $m^O=C$ in (\ref{eq:conservative_proj_not_consistent}) for a spatially constant field $C$, 
    \begin{equation}
        \intOmega{ \psi \rho^M m^M} = \intOmega{ \psi \rho^O C}, \quad \forall \psi \in \mathbb{V}_m^M
    \end{equation} 
    
    Given that $\rho^M$ is the Galerkin projection of $\rho^O$ (\ref{eq:galerkin_proj}) and $\mathbb{V}_m^M=\mathbb{V}_{\rho}^M$, we have
    \begin{align}
        \intOmega{ \psi \rho^M m^M} &= \intOmega{ \psi \rho^M C}, \quad \forall \psi \in \mathbb{V}_m^M, \\
        \rightarrow \intOmega{\psi \rho^M (m^M - C)} &= 0, \quad \forall \psi \in \mathbb{V}_m^M.
    \end{align}
    Choosing $\psi=m^M - C$ requires that $m^M = C$, given that $\rho^M>0$ (\ref{eq:positive_rho}). 
\end{pf}

We also introduce a conservative Galerkin projection that is consistent when $\rho^M$ and $m^M$ lie in different function spaces; this is required for projecting $m$ back to $\mathbb{V}_{\theta}$ after transport.

\begin{definition}[Consistent, conservative Galerkin projection]
The consistent, conservative Galerkin projection of $\overline{\mathcal{P}}_{\rho m}: \mathbb{V}_{m}^O, \mathbb{V}_{\rho}^O, \mathbb{V}_{\rho}^M \rightarrow \mathbb{V}_{m}^M$ computes $m^M\in\mathbb{V}_m^M$ from 
    \begin{equation} 
\intOmega{\psi \rho^M (m^M - m_{\overline{\Omega}})}= \intOmega{\psi \rho^O (m^O - m_{\overline{\Omega}} )}, \quad \forall \psi \in \mathbb{V}_m^M,
\label{eq:conservative_consistent_proj}
\end{equation}

\noindent where the constant
\begin{equation}
    m_{\overline{\Omega}} = \frac{\intOmega{m^O}}{\intOmega{}}
\label{eq:global_mean_m}
\end{equation}

\noindent is the global mean value of the mixing ratio.

\end{definition}

\begin{lemma}
\label{lemma:conservative_consistent_projection}
    The consistent, conservative Galerkin projection operator is conservative and consistent, given that $\rho^M$ is the Galerkin projection of $\rho^O$.
\end{lemma}

\begin{pf}
    For conservation, choose $\psi =1$ in (\ref{eq:conservative_consistent_proj}),
    \begin{equation}
        \intOmega{\rho^M m^M} = \intOmega{\rho^O m^O} + m_{\overline{\Omega}} \intOmega{(\rho^M-\rho^O)}.
    \end{equation}
    As $\rho^M$ is the Galerkin projection of $\rho^O$, $\intOmega{(\rho^M-\rho^O)}=0$ and we have mass conservation. \par
    For consistency, consider $m^O = C$ for some spatially constant $C$, which implies that $m_{\overline{\Omega}} = C$ (\ref{eq:global_mean_m}). Hence,
    \begin{equation} 
\intOmega{\psi \rho^M (m^M - C)}=0, \quad \forall \psi \in \mathbb{V}_{m}^M.
\end{equation}

\noindent Choosing $\psi = m^M - C$ requires that $m^M = C$, given that $\rho^M > 0$ (\ref{eq:positive_rho}). 
\end{pf}

The other component of the embedded DG method is the injection operator, which is used to map $m$ into the transport space. We first introduce the standard injection operator of \textcite{cotter2016embedded}, before describing a conservative version.

\begin{definition}[Injection]
    The injection operator, $\mathcal{I}: \mathbb{V}_{\theta} \rightarrow \hat{\mathbb{V}}_{\theta}$ identifies a fully discontinuous representation of a field from a partially continuous function space.
\end{definition}

As $\mathbb{V}_{\theta} \subset \hat{\mathbb{V}}_{\theta}$, the injection operator simply identifies a field in a larger function space and does not require the solution of a variational problem. The injection is sufficient for the advective scheme, but does not always ensure conservation of tracer density for the conservative scheme. To address this, we introduce a conservative injection operator which uses Galerkin projection.

\begin{definition}[Conservative injection]
    The conservative injection operator, $\mathcal{I}_{\rho m}: \mathbb{V}_{\theta}, \mathbb{V}_{\rho}, \hat{\mathbb{V}}_{\theta} \rightarrow \hat{\mathbb{V}}_{\theta}$ identifies a fully discontinuous representation of a mixing ratio ($m^M\in\hat{\mathbb{V}}_{\theta}$) from a partially continuous representation ($m^O\in\mathbb{V}_{\theta}$) by solving a conservative Galerkin projection (\ref{eq:conservative_proj_not_consistent}) with $\rho^O \in \mathbb{V}_{\rho}, \rho^M \in \hat{\mathbb{V}}_{\theta}$ and $\mathbb{V}_m^M=\hat{\mathbb{V}}_{\theta}$. As the function spaces of $\rho^M$ and $m^M$ are the same, the conservative injection operator is consistent by Lemma \ref{lemma:conservative_projection}. \par
\end{definition}

We now show that the conservative version of the vertically staggered scheme, which is summarised in Table \ref{table:vert_stag_operations}, is conservative and consistent.

\begin{theorem}[Vertically staggered scheme]
    The conservative vertically staggered scheme, which consists of (i) projection/injection from the original function spaces into the shared transport space, (ii) a full SSPRK3 timestep, and (iii) projection back to the original spaces, is conservative and consistent.
\end{theorem}

\begin{pf}
    By Lemmas \ref{lemma:conservative_projection} and \ref{lemma:conservative_consistent_projection}, steps (i) and (iii) are conservative and consistent. By Corollary \ref{corollary:ssprk_conservative_consistent}, each SSPRK3 stage, and thus the full timestep in step (ii), is conservative and consistent. 
\end{pf}

We note that we have global conservation for the vertically staggered scheme, but not local conservation; this is due to the vertical continuity of the $\mathbb{V}_{\theta}$ space and the consistent projection (\ref{eq:conservative_consistent_proj}) involving a global mean value. Local conservation could be obtained by using a vertically shifted mesh like in \textcite{bendall2023trilemma}.

\begin{table}[htpb]
\begin{center}
\caption[]{Sequence of operations for vertically staggered spaces, for the advective and conservative schemes, starting from original fields $\rho^n \in \mathbb{V}_{\rho}, ~m^n \in \mathbb{V}_{\theta}$ at time index $n$. The advective scheme solves (\ref{eq:rho_dg_upwind}) and (\ref{eq:advective_m_dg_upwind}), whilst the conservative scheme solves (\ref{eq:rho_dg_upwind}) and (\ref{eq:conservative_m_dg_upwind}). Fields with hats lie in the fully discontinuous transport space of $\mathbb{V}_T=\mathbb{\hat{V}}_{\theta}$. A dash in the conservative transport column denotes that the operation/function spaces are the same as the corresponding entry for the advective scheme.}
\label{table:vert_stag_operations}
\begin{tabular}{ |p{0.2\linewidth}||p{0.21\linewidth}|p{0.1\linewidth}||p{0.23\linewidth}|p{0.1\linewidth}| } 
\hline
Step & Advective scheme & Spaces & Conservative scheme & Spaces \\ 
\hline\hline
1. Project $\rho$ & $\hat{\rho}^n=\mathcal{P} [ \rho^n]$ & \mbox{$\mathbb{V}_{\rho} \rightarrow \mathbb{\hat{V}}_{\theta}$} & --- & --- \\ 
\hline
2. Inject or \newline conservatively inject $m$ & $\hat{m}^n = \mathcal{I} [ m^n]$ & \mbox{$\mathbb{V}_{\theta} \rightarrow \mathbb{\hat{V}}_{\theta}$} & $\hat{m}^n = \mathcal{I}_{\rho m } [ m^n, \rho^n, \hat{\rho}^n ] $ & $\mathbb{V}_{\theta},\mathbb{V}_{\rho},\mathbb{\hat{V}}_{\theta} \newline \rightarrow \mathbb{\hat{V}}_{\theta}$ \\ 
\hline
3. Transport & $[\hat{\rho}^{n+1},\hat{m}^{n+1}]$ = $\mathcal{T} [\hat{\rho}^{n},\hat{m}^{n}]$ & \mbox{$\mathbb{\hat{V}}_{\theta} \newline \rightarrow \mathbb{\hat{V}}_{\theta}$} & $[\hat{\rho}^{n+1},\hat{m}^{n+1}]$ = $\mathcal{T} [\hat{\rho}^{n}, \hat{\rho}^n \hat{m}^{n}]$ & --- \\
\hline
4. Project $\rho$ &\mbox{$\rho^{n+1}=\mathcal{P}[\hat{\rho}^{n+1}]$}& \mbox{$\mathbb{\hat{V}}_{\theta} \rightarrow \mathbb{V}_{\rho} $} & --- & --- \\
\hline
5. Project or \newline conservatively project $m$ & \mbox{$m^{n+1} = \mathcal{P} [\hat{m}^{n+1}] $} & \mbox{$\mathbb{\hat{V}}_{\theta} \rightarrow \mathbb{V}_{\theta}$} & $m^{n+1} = \overline{\mathcal{P}}_{\rho m} [\hat{m}^{n+1}, \hat{\rho}^{n+1}, \rho^{n+1}]$ & $\mathbb{\hat{V}}_{\theta},\mathbb{\hat{V}}_{\theta},\mathbb{V}_{\rho} \newline \rightarrow \mathbb{V}_{\theta}$ \\
\hline
\end{tabular}
\end{center}
\end{table}

\subsection{Monotonicity for co-located spaces}
\label{subsec:monotonicity}
The last key property of our conservative transport scheme is monotonicity. Due to the complexity of limiting partially continuous fields in a conservative manner, this work will only consider monotonicity in the co-located case. One approach to obtaining monotonicity in the vertically staggered case is the introduction of a vertically shifted mesh for limiting, similar to \textcite{bendall2023trilemma}. \par
Monotonicity of a field after DG transport can be enforced through a range of limiters \parencite[e.g.][]{kuzmin2010vertex, zhang2010maximum, cockburn2001runge}. An overview of DG limiting for atmospheric applications can be found in \textcite{nair2011emerging}. However, standard DG limiters applied to $m$ are designed to conserve $\int_e m \dx{V}$ in each element $e$, but do not enforce local tracer mass conservation, i.e. the element-wise preservation of $\int_e \rho m \dx{V}$. Limiters that conserve tracer mass in atmospheric contexts have been discussed in \textcite{schar1996synchronous} for flux-corrected transport, \textcite{guba2014optimization} for the spectral element method, and \textcite{bradley2019communication} for semi-Lagrangian transport. Our limiter has similarities to \textcite{guba2014optimization}, with the key differences being our DG application compared to continuous Galerkin finite elements, and that we explicitly blend a lower-order solution rather than selecting a limiter via optimisation. We note a similar conservative limiting example in \textcite{kuzmin2010failsafe} for the application of gas dynamics. \par
For atmospheric modelling, the task of conservative limiting also appears in physics--dynamics coupling applications when the physics and dynamics components are solved on meshes of differing resolution \parencite{herrington2019exploring,brown2024physics}. Our approach is similar to \textcite{brown2024physics}, which enforced non-negativity by blending a high-resolution, and possibly negative, field with a coarser resolution non-negative field. Instead, we consider the stricter case of monotonicity, and use a lowest-order mean mixing ratio (MMR) that is guaranteed to be monotonic under a CFL-type constraint. The MMR, $\overline{m} \in \mathbb{V}_{\rho,0}$, represents the mean value of $m\in \mathbb{V}_{\rho}$ in each element. Critically for mass conservation, $\overline{m}$ and $m$ represent the same tracer mass per cell, and so the convex combination of these fields is mass conserving.  \par

\begin{definition}[Mean mixing ratio]
    The mean mixing ratio (MMR), $\overline{m} \in \mathbb{V}_{\rho,0}$, of mixing ratio $m\in \mathbb{V}_{\rho}$ is defined in element $e$ as
    \begin{equation}
        \overline{m}_e = \frac{\int_e \rho m \dx{V}}{\int_e \rho \dx{V}}. \label{eq:mmr}
    \end{equation}   
\end{definition}

We now describe the bounds on the MMR that ensure monotonicity. Throughout, $\partial e$ denotes the boundary of a single cell $e$, with $\widehat{\bm{n}}$ the outwardly facing normal on $\partial e$. We then split the boundary into an inflow part $\partial e_{\mathrm{in}}$ containing the edges where $\bm{u\cdot}\widehat{\bm{n}}<0$, and an outflow part $\partial e_{\mathrm{out}}$ containing the edges where $\bm{u\cdot}\widehat{\bm{n}}\geq0$. With DG upwinding in the transport equations, values on $\partial e_{\mathrm{in}}$ are upwind traces from the neighbouring cells, whilst $\partial e_{\mathrm{out}}$ uses values from element $e$. Using $e\cup\partial e$ to denote $e$ along with the upwind traces on $\partial e_{\mathrm{in}}$, we define local extrema from the previous stage (pre-transported) mixing ratio,
\begin{equation}
    m^{i-1,\min}_{e\cup\partial e}=\min\left(m^{i-1,\min}_e,m^{i-1,\min}_{\partial e_{\text{in}}}\right), \quad m^{i-1,\max}_{e\cup\partial e}=\max\left(m^{i-1,\max}_e,m^{i-1,\max}_{\partial e_{\text{in}}}\right). 
\label{eq:local_monotonic_bounds}
\end{equation}

\noindent The extrema in each element are evaluated using the cell vertices, rather than the $\mathbb{V}_{\rho}$ nodes (Figure \ref{fig:function_spaces}), as the bilinear function in each element is guaranteed to have its maxima and minima at these locations.\par
As each SSPRK3 stage combines the time level $n$ field and a forward Euler update (\ref{eq:ssprk3}), each stage has local lower and upper bounds for monotonicity of
\begin{equation}
    m^{i,-}_e = \min \bracked{m^{i-1,\min}_{e\cup\partial e}, m^{n, \min}_e}, \quad m^{i,+}_e = \max \bracked{m^{i-1,\max}_{e\cup\partial e}, m^{n, \max}_e}.
\label{eq:high_low_bounds}
\end{equation}

\noindent Monotonicity of the MMR is enforced by requiring that
\begin{equation}
    m^{i,-}_e \leq \overline{m}^{i}_e \leq m^{i,+}_e,
\label{eq:mmr_bounds}
\end{equation}

\noindent which is obeyed subject to a CFL-type condition on the outflow in each element. This bound is proved in Theorem \ref{theorem:monontonic_mmr} and Corollary \ref{corollary:stage_mmr_monotonic} in Appendix \ref{append:monotonic_mmr}. \par
We now describe the blended field that is obtained after applying the MMR limiter.

\begin{definition}[Blended mixing ratio]
    The blended mixing ratio $m^*_e$ in element $e$ is
    \begin{equation}
		m^{*}_e = (1-\Lambda_e)m^{i}_e + \Lambda_e\overline{m}^{i}_e,
    \label{eq:blended_m}
	\end{equation}
    \noindent where $m^{i}_e$ is the mixing ratio after stage $i$ of the SSPRK3 scheme, and $\Lambda_e \in [0,1]$ is the blending weight for element $e$. The value of $\Lambda_e$ which enforces monotonicity is
\begin{equation}
    \Lambda_e = \max\left(\Lambda^{\min}_e, \Lambda^{\max}_e \right),
\label{eq:lambda_e}
\end{equation}
where
\begin{equation}
    \Lambda^{\min}_e = \begin{dcases}
         \frac{m^{i,-}_e - m^{i,\min}_e}{\overline{m}_e - m^{i,\min}_{e}}, & \text{if} ~m^{i,\min}_e < m^{i,-}_e, \\
         0, & \text{otherwise},
    \end{dcases}
    \qquad
    \Lambda^{\max}_e = \begin{dcases}
         \frac{m^{i,\max}_{e} - m^{i,+}_e }{m^{i,\max}_{e} - \overline{m}_e}, & \text{if} ~m^{i,\max}_e > m^{i,+}_e, \\
         0, & \text{otherwise}.
    \end{dcases}
\label{eq:lambda_min}
\end{equation}
\end{definition}

We note that the blending operation (\ref{eq:blended_m}) follows the general structure of \textcite{barth1989design}. $\Lambda^{\min}_e$ enforces the lower bound by ensuring that, if the mixing ratio falls below the lower bound in element $e$, then $m^{i,\min}_e$ is increased to $m^{i,-}_e$. This leads to a blending weight of
\begin{subequations}
    \begin{align}
        m^{i,-}_e & = (1-\Lambda^{\min}_e) m^{i,\min}_e + \Lambda^{\min}_e \overline{m}_e, \\
        \rightarrow \Lambda^{\min}_e & = \frac{m^{i,-}_e - m^{i,\min}_e}{\overline{m}_e - m^{i,\min}_{e}}.
    \end{align}
\end{subequations}

\noindent Similarly, $\Lambda^{\max}_e$ enforces the upper bound by ensuring that, if the mixing ratio field exceeds the upper bound, then $m^{i,\max}_e$ is reduced to $m^{i,+}_e$ after blending. \par
Given the MMR bounds for monotonicity (\ref{eq:mmr_bounds}), $\Lambda^{\min}_e, \Lambda^{\max}_e \in [0,1]$ in the absence of numerical error. The denominator is never evaluated as zero, as for example in $\Lambda^{\min}_e$, $\overline{m}_e - m^{i,\min}_{e} = 0$ would require that $\overline{m}_e = m^{i,\min}_{e} = m^{i,\max}_{e}$, and $\overline{m}_e \geq m_e^{i,-}$ (\ref{eq:mmr_bounds}), but in this case, $\Lambda^{\min}_e=0$ is enforced by the conditional statement (\ref{eq:lambda_min}). However, to avoid issues with numerical error, the $[0,1]$ bounds are strictly enforced in the numerical implementation. \par
We now summarise the MMR limiter process at each stage:
\begin{itemize}
    \item Calculate the pre-transport local extrema ($m^{i-1,\min}_{e\cup\partial e}$, $m^{i-1,\max}_{e\cup\partial e}$) (\ref{eq:local_monotonic_bounds}), then the stage $i$ local bounds ($m^{i,-}, m^{i,+}$) (\ref{eq:high_low_bounds})
    \item Solve the SSPRK3 stage to obtain the pre-limited field $m^i$
    \item Compute the MMR $\overline{m}^i$ (\ref{eq:mmr})
    \item Compute the blending weights $\Lambda_e \in [0,1]$ (\ref{eq:lambda_e})
    \item Compute the blended, monotonic, field $m^*$ (\ref{eq:blended_m}) 
\end{itemize}

The MMR limiter can be used in systems with multiple mixing ratios. In this case, it is important to preserve any linear correlations, where a combination of mixing ratio fields forms a constant field that should be conserved \parencite{lauritzen2012evaluating}. Preservation of linear correlations is obtained by applying the same maximal blending function to all mixing ratios \parencite{brown2024physics}, and a proof of this is provided in Theorem \ref{theorem:maximal_blending_fun} in Appendix \ref{append:proof_blending_func}. The preservation of linear correlations is examined in the terminator toy test of Section \ref{subsec:terminator}.

\begin{definition}[Maximal Blending Function]
    Given individual blending functions $\Lambda_{k,e}$ for mixing ratio $k$ in element $e$, the maximal blending function $\Lambda^M_e$ in element $e$ is given by
\begin{equation}
    \Lambda^M_e = \max_k \{ \Lambda_{k,e}\}.
\label{eq:multiple_blend_func}
\end{equation}
\noindent Using the maximum ensures that monotonicity is enforced for all mixing ratios.
\end{definition}

Now that we have shown the key properties of our conservative transport scheme, we move to numerical tests, starting with the co-located configuration.

\section{Results for co-located spaces}
\label{section:colocate_spaces}
In the co-located configuration, the dry density and mixing ratio live in the same fully discontinuous function space, $\mathbb{V}_{\rho}$. The transport operator is applied in $\mathbb{V}_T=\mathbb{V}_{\rho}$, thus no mapping between function spaces is required. We provide three tests in the co-located setup. The first is a transport-only test to examine conservation and consistency. The second test checks that the MMR limiter enforces monotonicity. The third test includes physics terms and multiple mixing ratios, and examines conservation, monotonicity, and the preservation of linear correlations. All co-located tests use a spherical domain with an equiangular gnomonic cubed-sphere mesh \parencite{ronchi1996_equiangular}.

\subsection{Transport-only test}
\label{subsec:co_located_test}
We first apply a transport test on the sphere from \textcite{nair2010class}. The use of only transport terms allows a direct comparison of the conservative and advective form schemes for the tracer transport. We run two versions of the test: first, the convergence setting, which tests mass conservation and the order of the scheme convergence, and second, the consistency setting, which initialises a constant mixing ratio and checks whether this is preserved over time. \par
The velocity in this test is prescribed as the divergent flow from \textcite{lauritzen2012standard}, leaving the density and mixing ratio as the only prognostic variables. The flow field is defined such that a scalar field should be transported back to its original location at the end time. Hence, any deviation from the initial condition provides a measure of the transport error. In the convergence setting, two Gaussian perturbations are applied to the mixing ratio field, whilst the consistency test uses a constant mixing ratio field with a perturbed density. Further details on the test case setup are found in Appendix \ref{appendix:setup1}. \par
\begin{figure}
    \centering
    \includegraphics[width=0.9\linewidth]{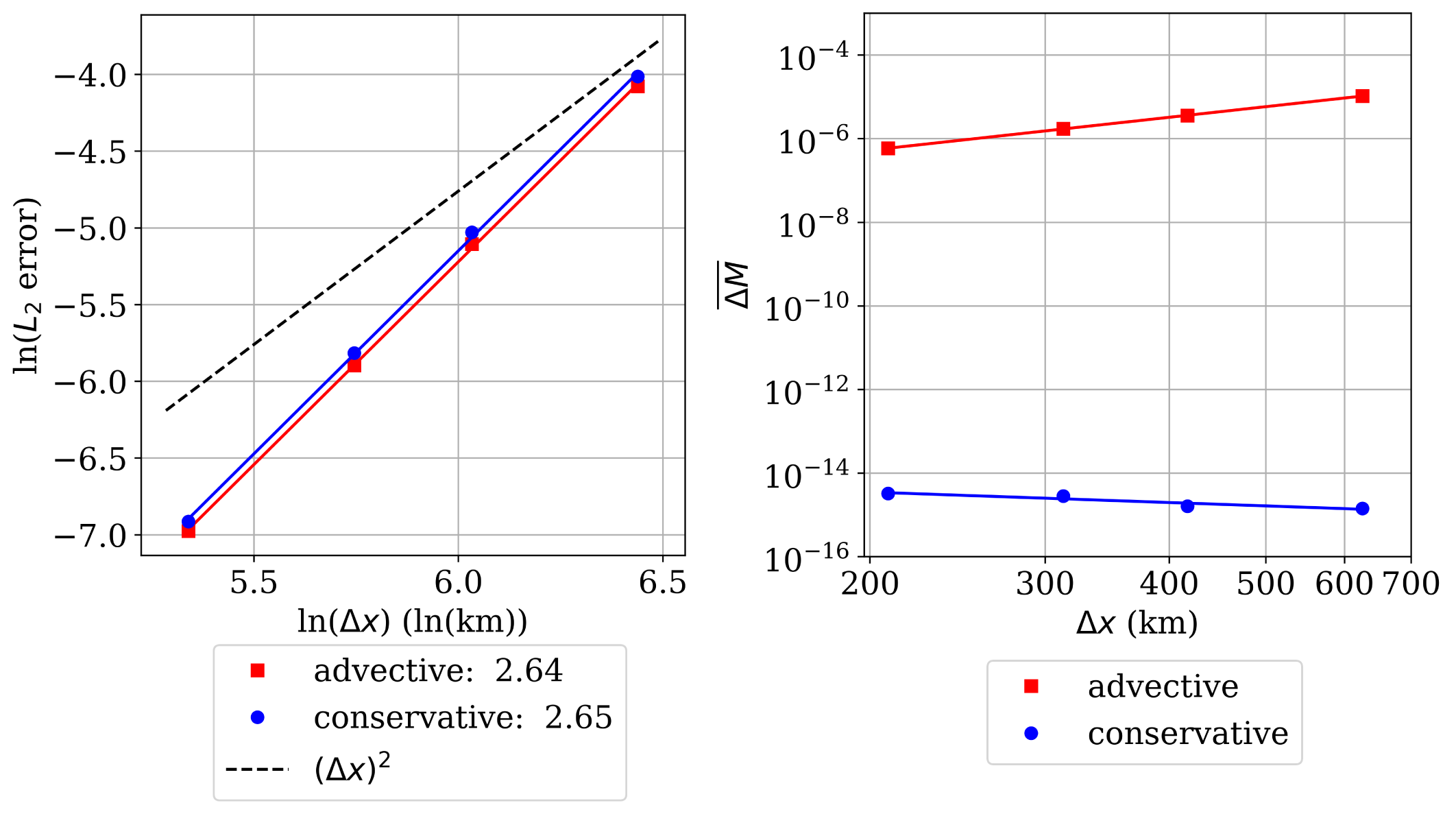}
    \caption{(Left) $L_2$ error in the mixing ratio (\ref{eq:m_L2_err}), and (right) mean mass error (\ref{eq:mean_td_error}) in the co-located transport-only tests. A characteristic cell length for the cubed-sphere is computed as $\Delta x= \pi R/(2 N_e)$, where $N_e$ is the number of cells per panel edge. Legend values in the $L_2$ error plot (left) denote the approximate order of convergence from a line of best fit. The advective and conservative schemes have virtually the same order of convergence, with the conservative scheme obeying mass conservation.}
    \label{fig:NL_sphere_error_td}
\end{figure}
\begin{figure}
    \centering
    \includegraphics[width=0.6\linewidth]{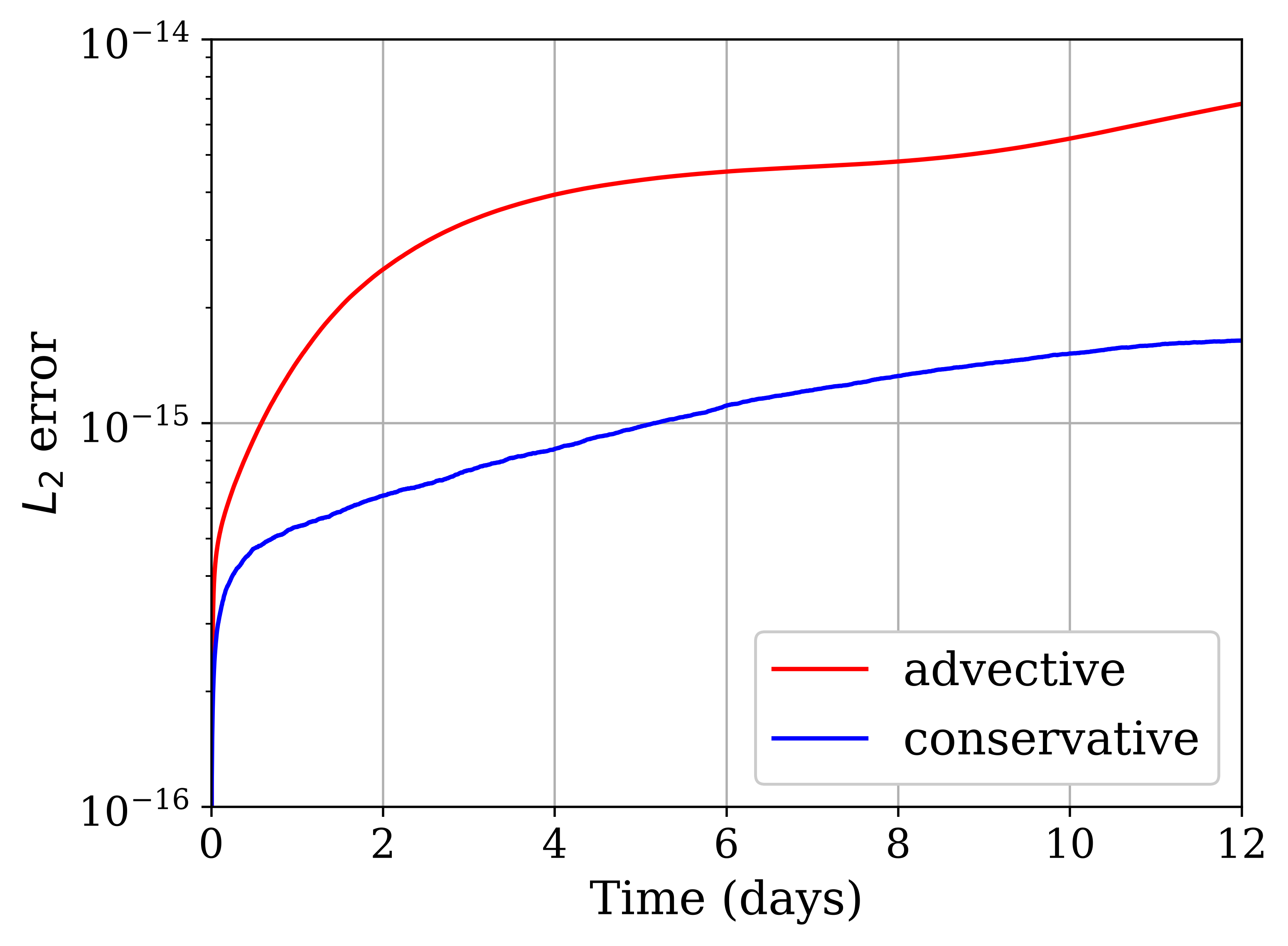}
    \caption{$L_2$ error in the mixing ratio (\ref{eq:m_L2_err}) over time in the co-located transport-only test with the consistency configuration, in which the mixing ratio is held constant while the density varies spatially. Both schemes respect the consistency property.}
    \label{fig:NL_sphere_consistent}
\end{figure}

For the convergence setting, Figure \ref{fig:NL_sphere_error_td} shows the $L_2$ error (\ref{eq:m_L2_err}) in the final mixing ratio fields, with both schemes showing the same convergence rates, which are better than second-order. Figure \ref{fig:NL_sphere_error_td} also shows the mean mass error (\ref{eq:mean_td_error}) over the simulation. The conservative tracer transport scheme preserves the tracer density to near machine precision, indicating mass conservation. Contrastingly, there are significant mass discrepancies with the advective scheme. For the consistency setting, Figure \ref{fig:NL_sphere_consistent} shows that both the advective and conservative tracer methods are consistent, as they preserve the constant mixing ratio to near machine precision. \par

\subsection{MMR limiter test}
\begin{figure}
    \centering
    \includegraphics[width=0.8\linewidth]{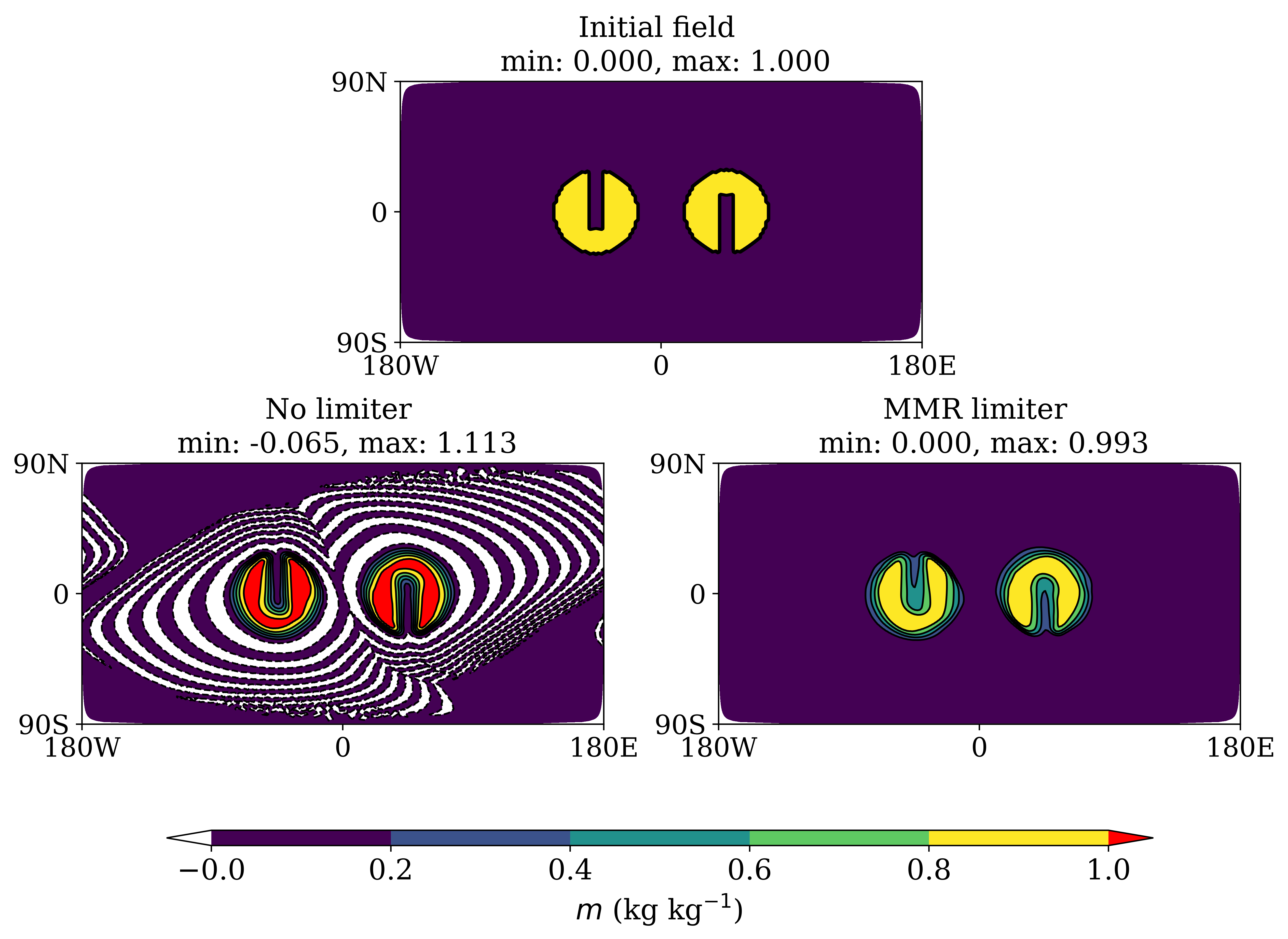}
    \caption{The initial mixing ratio field (top) from the MMR transport test, along with the final mixing ratio fields with no limiter (bottom left) and the MMR limiter (bottom right). The MMR limiter successfully enforces monotonicity.}
    \label{fig:mmr_test_m_fields}
\end{figure}

\begin{figure}
    \centering
    \includegraphics[width=0.7\linewidth]{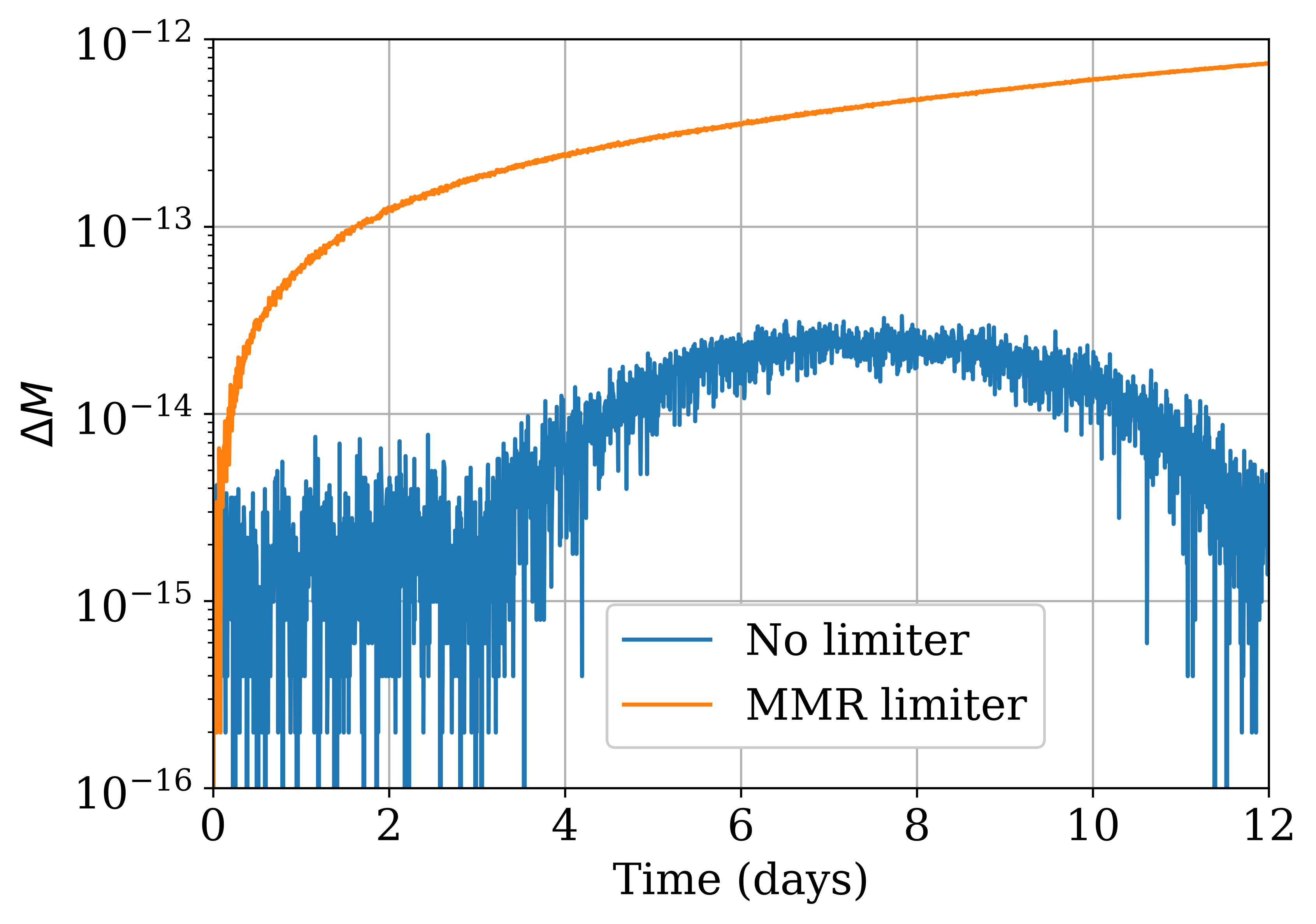}
    \caption{Relative change in tracer mass $\Delta M$ (\ref{eq:td_error_over_time}) over time in the MMR limiter test. The MMR limiter keeps $\Delta M$ below $10^{-12}$ during the test.}
    \label{fig:mmr_test_tracer_density}
\end{figure}

We now test the MMR limiter in a similar setup to the convergence setting of the transport-only test, but now with a slotted cylinder initial condition, as described in Appendix \ref{appendix:setup2}. Discontinuities in the initial mixing ratio field, which is zero outside of the two slotted cylinders, mean that under- and overshoots are easily generated during transport. The conservative tracer transport scheme is used for the mixing ratio, and simulations with no limiting and with the MMR limiter are compared. The unlimited case allows undershoots ($m<0$) and overshoots ($m>1$) in the mixing ratio, whereas the monotonic limiter enforces $m \in [0,1]$ throughout the domain (Figure \ref{fig:mmr_test_m_fields}). The MMR limiter achieves monotonicity whilst preserving the tracer mass to near machine precision (Figure \ref{fig:mmr_test_tracer_density}).

\clearpage

\subsection{Terminator toy test}
\label{subsec:terminator}
The last co-located test examines the preservation of a linear correlation between two mixing ratios. We apply the terminator toy chemistry system of \textcite{lauritzen2015terminator}, which considers the transport of two molecular species, specified through mixing ratios\footnote{Note that \textcite{lauritzen2015terminator} states that  $X$ and $X_2$ are number densities, with $X_T$ the total number of atoms of $X$ per molecule of dry air. However, we will treat these quantities as mixing ratios, as per the design of this test for dynamical cores; this means we still quote units of kg kg$^{-1}$.} of $X$ and $X_2$. $X$ is a single-atom species and $X_2$ is comprised of two atoms, so the total mixing ratio is $X_T = X + 2X_2$. The exchange of molecules between species $X$ and $X_2$ is governed by a pair of ODEs that preserve the total tracer mass of $M=\intOmega{\rho (X + 2X_2)}$. The addition of these chemical equations as forcing to the transport equations is considered the `physics' applied to a dynamical core. The presence of physics terms makes the equations more challenging to solve and provides a potential source of negative mixing ratios. Without a limiter, this case is likely to break monotonicity, due to steep gradients in the initial mixing ratio fields at the terminator. \par 
Consider co-located variables of $\rho, X, X_2 \in \mathbb{V}_{\rho}$. The advective scheme solves transport equations of
\begin{subequations}
        \begin{align}
         \pfrac{\rho}{t} + \bm{\nabla\cdot}\left(\rho \vec{u}\right) &= 0, \\
         \pfrac{X}{t} + (\vec{u} \cdot \bm{\nabla}) X &= 0,
        \\
         \pfrac{X_2}{t} + (\vec{u} \cdot \bm{\nabla}) X_2 &= 0,
        \end{align}
\label{eq:terminator_toy_adv}    
\end{subequations}

\noindent whilst the conservative tracer scheme solves

\begin{subequations}
        \begin{align}
         \pfrac{\rho}{t} + \bm{\nabla\cdot}\left(\rho \vec{u}\right) &= 0, \\
         \pfrac{(\rho X)}{t} + \bm{\nabla\cdot}\left(\rho X \vec{u}\right) &= 0, \label{eq:terminator_X_transport}
        \\
         \pfrac{(\rho X_2)}{t} + \bm{\nabla\cdot}\left(\rho X_2 \vec{u}\right) &= 0.
        \label{eq:terminator_X2_transport}
        \end{align}
\label{eq:terminator_toy}  
\end{subequations}

The physics terms in this system are the coupled chemical ODEs governing the molecular interactions,
\begin{subequations}
\begin{align}
        \der{X}{t} &= 2f_X, \label{eq:terminator_toy_ode_X}  \\
        \der{X_2}{t} &= -f_X, \label{eq:terminator_toy_ode_X2}
\end{align}
\label{eq:terminator_toy_ode}
\end{subequations}

\noindent where $f_X = k_1 X_2 - k_2 X^2$, and $k_1, k_2$ are reaction rates given in Appendix \ref{appendix:setup3}. Combining (\ref{eq:terminator_toy_ode_X}) with twice (\ref{eq:terminator_toy_ode_X2}) shows that $\mathrm{d}X_T/\mathrm{d}t = 0$, and hence, the chemistry ODEs are conservative. The spatial profile of $k_1$ contains strong gradients, making this a challenging test for transport schemes as numerical error can easily generate spurious under- and overshoots. As the chemical ODEs (\ref{eq:terminator_toy_ode}) are too stiff to solve explicitly with a reasonable timestep size, $f_X$ is computed from implicit forcing expressions given in Appendix D of \textcite{lauritzen2015terminator}. To ensure that the physics step does not introduce negative values in the mixing ratio fields, $f_X$ is evaluated at the cell vertices and is bounded as
\begin{equation}
    f_X = \begin{dcases}
        \max \bigg{\{} f_X, -\frac{X}{2 \Delta t} \bigg{\}}, ~\text{if} ~f_X<0, \\
        \min \bigg{\{} f_X, \frac{X_2}{\Delta t} \bigg{\}}, ~\text{else} . 
    \end{dcases}
\label{eq:fX_bounds}
\end{equation}

\noindent Updates to $X$ and $X_2$ from the physics are performed after each full SSPRK3 transport step, with the bounding of $f_X$ (\ref{eq:fX_bounds}) ensuring non-negative mixing ratios at the start of each SSPRK3 step. For other physics schemes that might generate negative values, any negativity in the mixing ratios should be removed using a mass-conserving limiter before beginning the transport timestep. \par
The initial condition for the mixing ratios is such that the total tracer mixing ratio is spatially constant: $X_T(0)=X(0) + 2X_2(0)=4\times 10^{-6} ~\text{kg} ~\text{kg}^{-1}$. Hence, a scheme that preserves linear correlations will ensure that $X_T(t) = X_T(0)$. To preserve linear correlations with the MMR, the maximal blending field $\Lambda^M$ (\ref{eq:multiple_blend_func}) is used to limit $X$ and $X_2$ at each SSPRK3 stage. For the advective scheme, the vertex-based slope limiter of \textcite{kuzmin2010vertex} is applied to the two mixing ratios. The density is unlimited in both schemes. Further details about the test case setup are found in Appendix \ref{appendix:setup3}. \par
\begin{figure}[htpb]
        \centering
         \includegraphics[width=0.6\textwidth]{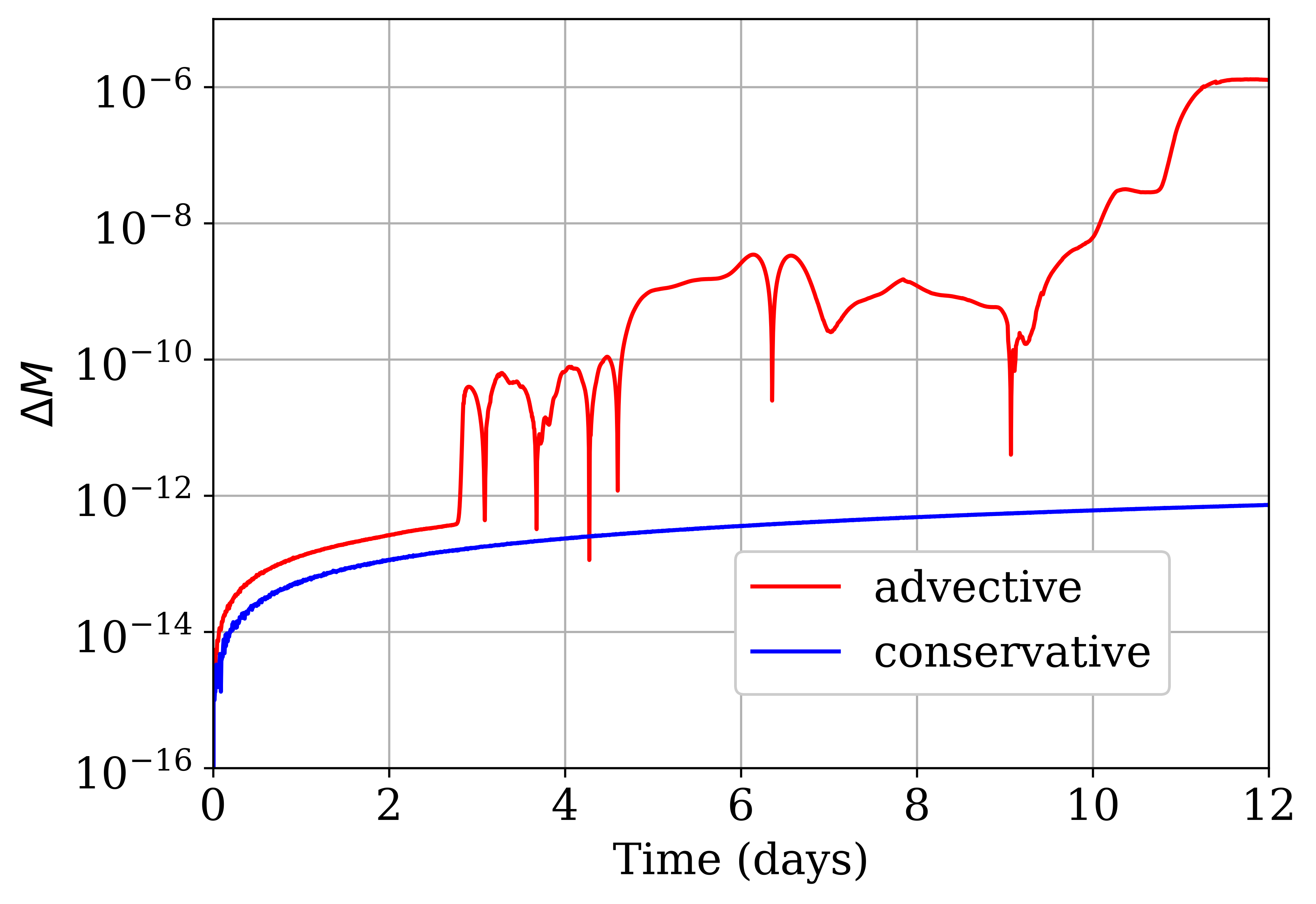}
        \caption{Relative change in total tracer mass $\Delta M$ (\ref{eq:td_error_over_time}), where $M=\intOmega{\rho(X+2X_2)}$, for the advective and conservative tracer transport schemes during the terminator toy test. The conservative scheme preserves the tracer mass better than the advective scheme. }
    \label{fig:terminator_toy_tracer_density}
\end{figure}

Figure \ref{fig:terminator_toy_fields} shows the initial condition for $\rho, X$ and $X_2$ in the terminator toy test, along with the final states for the advective and conservative schemes. The final  fields from both schemes are visually indistinguishable, and both have monotonic final states for the mixing ratios, although the advective scheme is not guaranteed to maintain monotonicity during the simulation. There is a difference in the mixing ratio extrema, with the advective scheme generating a smaller maximum value of $X_2$, and a larger minimum value of $X$ in the final states, which is a reflection of the advective scheme limiter being more diffusive than the MMR limiter in this test. The advective scheme does not preserve the linear correlation in the mixing ratios, as shown by the nonzero $X_T-X_T(0)$ field at the end of the simulation (bottom row of Figure \ref{fig:terminator_toy_fields}). Conversely, the magnitude of $X_T-X_T(0)$ with the conservative scheme is at machine precision, meaning that it has preserved the constant total mixing ratio. Mass conservation of the conservative transport scheme is also shown through a time series of relative change in tracer mass in Figure \ref{fig:terminator_toy_tracer_density}. Both methods have similar mass conservation up to the end of day two, but beyond this the conservative scheme maintains total tracer mass much better than the advective scheme. \par

\begin{figure}[htpb]
        \centering
         \includegraphics[width=\textwidth]{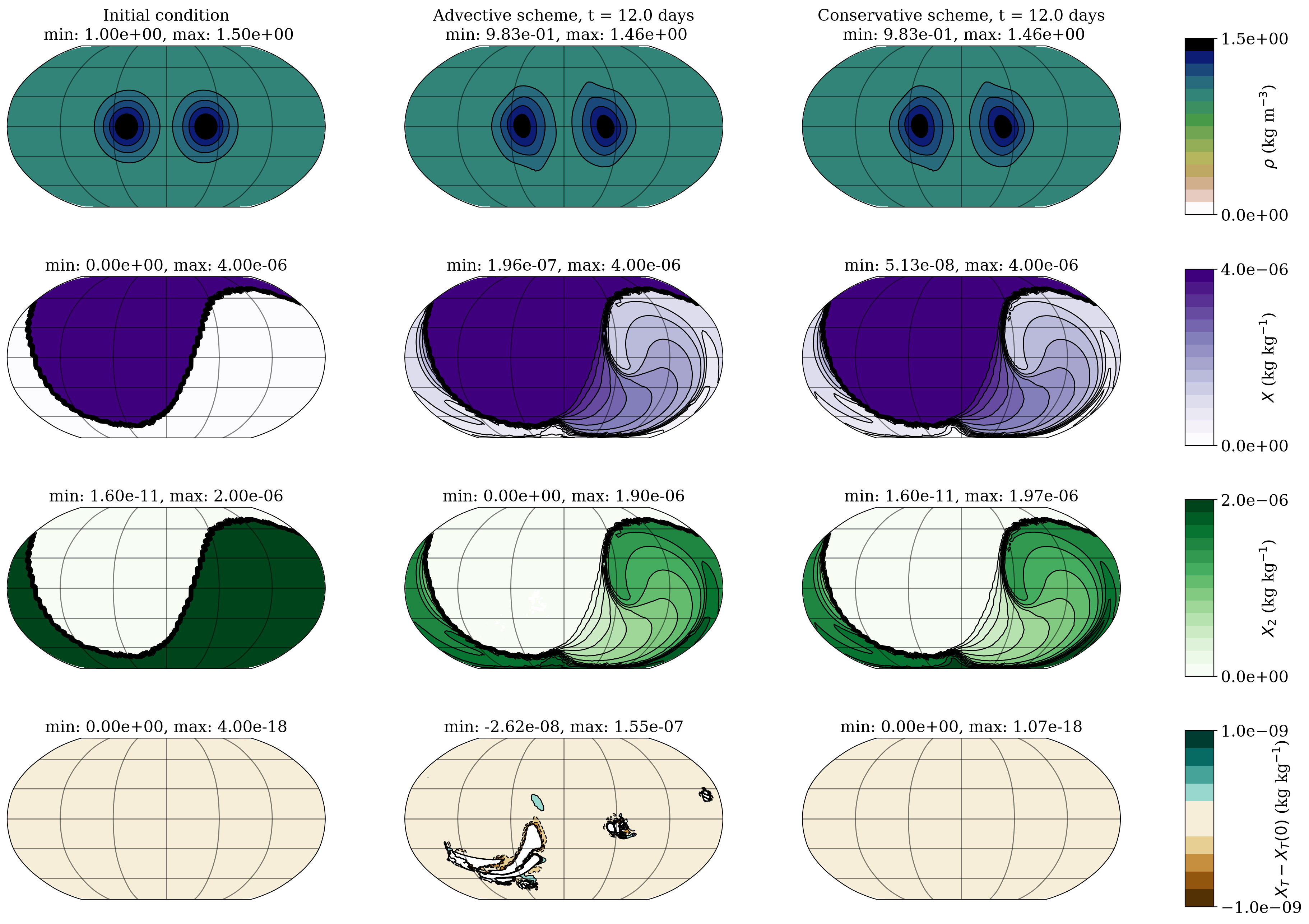}
        \caption{Plots of the density and mixing ratios for the tracer species $X$ and $X_2$ in the terminator toy test, visualised with a Robinson projection of the sphere. The left-hand column shows the initial condition, the middle column the final state with the advective scheme, and the right-hand column the final state with the conservative tracer scheme. The first row shows the density fields, the middle two rows plot $X$ and $X_2$, and the last row plots the change in the total mixing ratio, $X_T-X_T(0)$, where $X_T(0) = 4\times 10^{-6} ~\text{kg} ~\text{kg}^{-1}$ is a constant field. The conservative scheme preserves the constant total mixing ratio.}
    \label{fig:terminator_toy_fields}
\end{figure}

\clearpage

\section{Results for vertically staggered spaces}
\label{section:vert_stag_spaces}
We now test the conservative transport scheme when the density and mixing ratio lie in different function spaces: $\rho \in \mathbb{V}_{\rho}, m \in \mathbb{V}_{\theta}$. This requires injection and projection operators to map the mixing ratio and density to a shared transport space and back to their original spaces (Table \ref{table:vert_stag_operations}). These tests use a vertical slice domain, with the thermodynamic function space $\mathbb{V}_{\theta}$ being discontinuous in the horizontal dimension and continuous in the vertical (Figure \ref{fig:function_spaces}). \par

\subsection{Transport-only test}
\label{subsec:vert_staggered_test}
The first test with vertically staggered spaces uses a vertical slice modification of the co-located transport-only test (Section \ref{subsec:co_located_test}), as proposed in \textcite{bendall2023trilemma}. Like in the co-located transport-only tests, we run two configurations --- the convergence setting, with a perturbed mixing ratio field to test mass conservation, and the consistency setting with a constant mixing ratio that should remain constant. Further details of the test case setup can be found in Appendix \ref{appendix:setup4}. \par
\begin{figure}[htpb]
    \centering
    \includegraphics[width=0.9\linewidth]{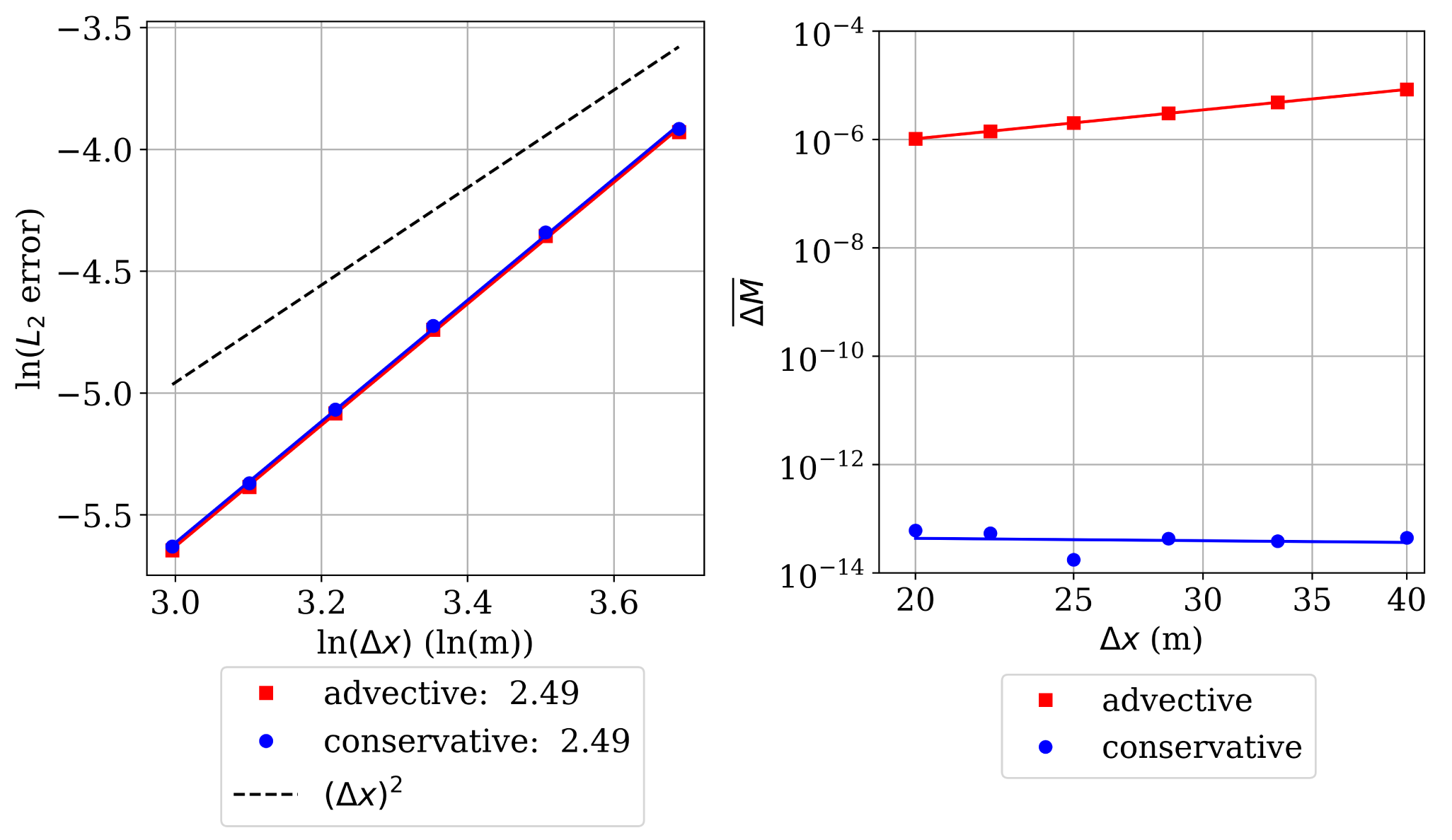}
    \caption{$L_2$ error in the mixing ratio (\ref{eq:m_L2_err}) (left) and mean mass error (\ref{eq:mean_td_error}) (right) in the vertically staggered transport-only test. Legend values in the $L_2$ error plots denote the approximate order of convergence from a line of best fit. Both schemes have the same rate of convergence, but the conservative scheme maintains much better mass conservation.}
    \label{fig:NL_slice_mx_err_and_td}
\end{figure}
\begin{figure}[htpb]
    \centering
    \includegraphics[width=0.6\linewidth]{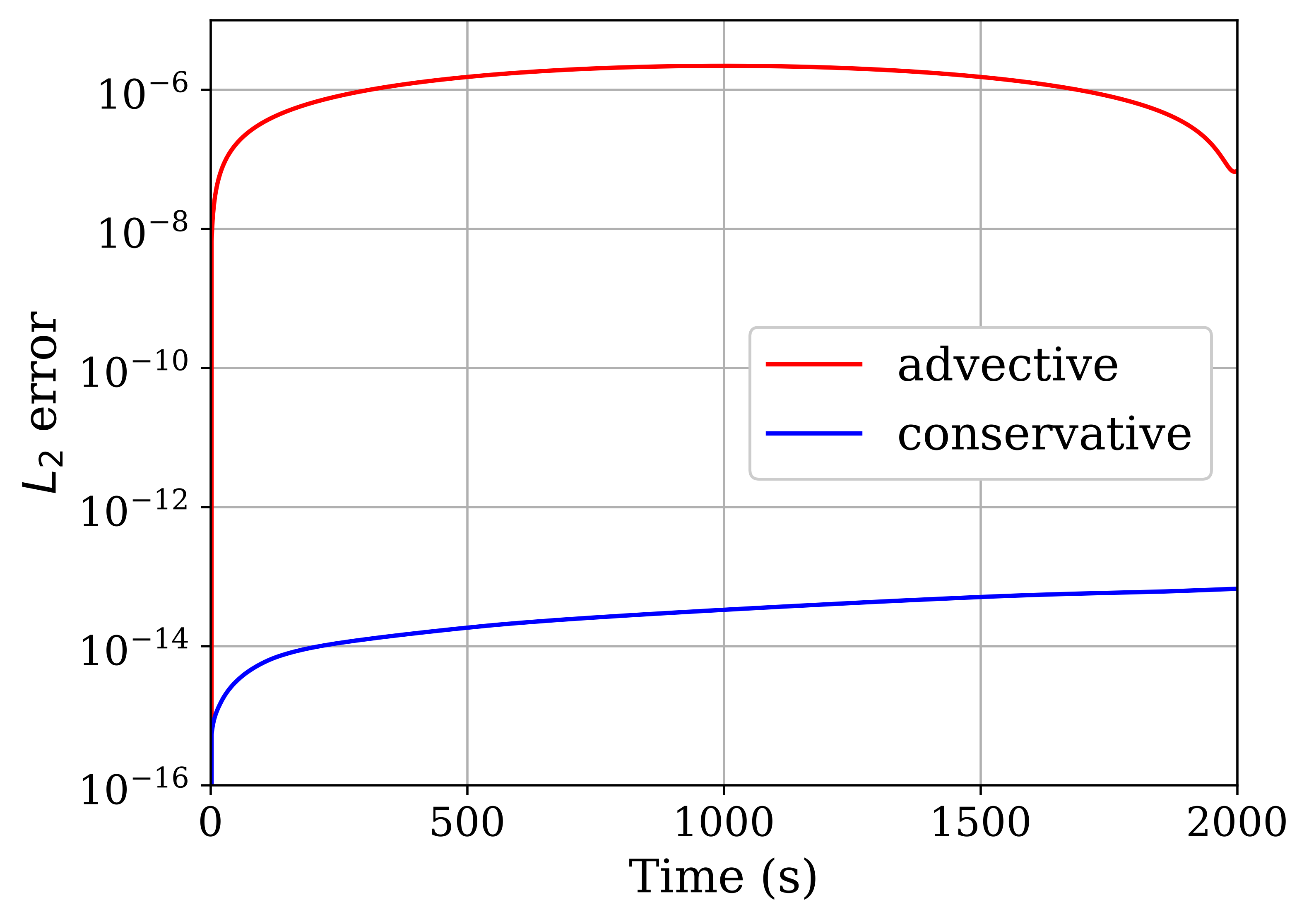}
    \caption{$L_2$ error in the mixing ratio (\ref{eq:m_L2_err}) in the consistency setting of the vertically staggered transport-only test, in which the mixing ratio is held constant while the density is spatially varying. The conservative scheme is much more consistent than the advective scheme.}
    \label{fig:NL_slice_consistent}
\end{figure}

For the convergence test, Figure \ref{fig:NL_slice_mx_err_and_td} shows the $L_2$ error in the final mixing ratio field and the mean mass error. As with the co-located test, both transport schemes have the same, above second-order, convergence rate. The conservative tracer scheme has mass conservation that is close to machine precision, showing that the new scheme can conserve mass when there is a vertical staggering of $\rho$ and $m$. For the consistency test, Figure \ref{fig:NL_slice_consistent} shows that the conservative transport scheme preserves the constant mixing ratio field much better than the advective scheme. The advective scheme is much less consistent than it was with co-located function spaces (Figure \ref{fig:NL_sphere_consistent}). \par

\clearpage

\subsection{Dynamical core test}
\label{sec:dycore_test}
As a final test, we apply the conservative tracer transport scheme in a nonhydrostatic dynamical core setup. This solves the compressible Euler equations in a vertical slice domain, without the Coriolis force, and coupled to a single physics scheme of saturation adjustment. We follow the governing equations outlined in \textcite{bendall2020compatible}, but neglect the rain moisture species. \par
The prognostic variable set is $[\vec{u},\rho,\theta_{\text{vd}},m_{\text{v}},m_{\text{c}}]^T$, with $\vec{u}=[u,w]^T$ the velocity in the vertical slice, $\rho$ the dry air density, $\theta_{\text{vd}}$ the virtual dry potential temperature, and $m_{\text{v}}$ and $m_{\text{c}}$ the mixing ratios for the moisture species of water vapour and cloud water, respectively. To conserve water mass, we require that $M = \intOmega{\rho(m_{\text{v}} + m_{\text{c}})}$ is constant in time. With a vertical staggering, we have $\theta_{\text{vd}}, m_{\text{v}}, m_{\text{c}} \in \mathbb{V}_{\theta}$ and $\rho \in \mathbb{V}_{\rho}$. Unlike the previous transport tests where the velocity was prescribed, $\vec{u}$ is now a prognostic field defined in a Raviart-Thomas function space.  \par
The governing equations, using the advective form for the transport of mixing ratios, are
\begin{subequations}
\begin{align}
    \pfrac{\vec{u}}{t} + (\bm{\nabla} \times \vec{u}) \times \vec{u} + \frac{1}{2} \bm{\nabla} |\vec{u}|^2  + \frac{c_{\text{p}} \theta_{\text{vd}}}{1 + m_{\text{v}} + m_{\text{c}}} \bm{\nabla} \Pi + g \vec{\hat{k}} &= \vec{0}, \label{eq:euler_eqs_u} \\
    \pfrac{\rho}{t} + \bm{\nabla} \cdot (\vec{u} \rho) &= 0, \\
    \pfrac{\theta_{\text{vd}}}{t} + (\vec{u} \cdot \bm{\nabla}) \theta_{\text{vd}} + \mu \theta_{\text{vd}} (\bm{\nabla} \cdot \vec{u}) &= S_{\theta}, \\
    \pfrac{m_{\text{v}}}{t} + (\vec{u} \cdot \bm{\nabla}) m_{\text{v}} &= S_{\text{v}}, \label{eq:w_v_eq} \\
    \pfrac{m_{\text{c}}}{t} + (\vec{u} \cdot \bm{\nabla}) m_{\text{c}} &= S_{\text{c}},  \label{eq:c_w_eq}
\end{align}
\label{eq:euler_eqs}
\end{subequations}

\noindent with $c_{\text{p}}$ the specific heat capacity of dry air, $\vec{\hat{k}}$ the unit vector in the vertical $z$-coordinate, and $\Pi$ the Exner pressure. Note the use of the vector-invariant form for momentum transport (\ref{eq:euler_eqs_u}). The coefficient $\mu$ in the thermodynamic equation is 
\begin{equation}
    \mu = \frac{R_{\text{m}}}{c_{\text{vml}}} - \frac{R_{\text{d}} c_{\text{pml}} }{ c_{\text{p}} c_{\text{vml}} },
\end{equation}

\noindent with $R_{\text{m}}$ a specific gas constant, $R_{\text{d}}$ the specific gas constant for dry air, and $c_{\text{pml}}, c_{\text{vml}}$ representing specific heat capacities for moist air, at constant pressure and constant volume, respectively. \par
The $S$ terms in (\ref{eq:euler_eqs}) comprise the saturation adjustment physics scheme, which computes the condensation of water vapour to cloud water ($S_{\text{c}}$), the evaporation of cloud water to water vapour ($S_{\text{v}}$), and the effect of latent heat change on the virtual dry potential temperature ($S_{\theta}$). The physics scheme ensures that mass is not lost during the water phase changes, i.e. $S_{\text{v}} + S_{\text{c}} = 0$. General details of the saturation adjustment scheme can be found in \textcite{bryan2002benchmark}, with the finite element formulation we use given in \textcite{bendall2020compatible}. \par 
The compressible Euler equation set is closed with an equation of state for the Exner pressure,
\begin{equation}
    \Pi = \bracfrac{\rho R_{\text{d}} \theta_{\text{vd}}}{p_0}^{\frac{\kappa}{1-\kappa}},
\label{eq:exner}
\end{equation}

\noindent where $p_0 = 1000 ~\text{hPa}$ is a reference pressure, and $\kappa = R_{\text{d}}/c_{\text{p}}$. The Exner pressure is computed diagnostically from the $\rho$ and $\theta_{\text{vd}}$ fields, which are vertically staggered on the Charney-Phillips grid. We choose $\Pi \in \mathbb{V}_{\rho}$ to co-locate the Exner pressure with the density, so the projection $\mathcal{P}:\mathbb{V}_{\theta} \rightarrow \mathbb{V}_{\rho}$ is applied to $\theta_{\text{vd}}$ before evaluating $\Pi$ (\ref{eq:exner}). \par
As the compressible Euler equations support fast acoustic waves, applying an explicit timestepper like SSPRK3 to the complete equation set would be restricted to much smaller timesteps for stability than in the previous tests. As such, we use a semi-implicit quasi-Newton (SIQN) timestepper \parencite{wood2014inherently,melvin2019cartesian}, which mirrors the timestepping methodology of GungHo. Pseudocode for SIQN can be found in \textcite{hartney2025exploring}. The SIQN timestepper consists of an outer loop for the slower transport terms, and an inner loop where the fast dynamics is solved implicitly. For the compressible Euler equations, the pressure gradient ($\bm{\nabla} \Pi$) and gravity ($g \vec{\hat{k}}$) terms from the momentum equation (\ref{eq:euler_eqs_u}) are solved in the inner loop. The outer loop solves a transport equation for each prognostic variable, with the conservative scheme using conservative transport equations for the moisture species,
\begin{subequations}
\begin{align}
    \pfrac{(\rho m_{\text{v}})}{t} + \bm{\nabla} \cdot (\rho m_{\text{v}} \vec{u}) &= 0, \\
    \pfrac{(\rho m_{\text{c}})}{t} + \bm{\nabla} \cdot (\rho m_{\text{c}} \vec{u}) &= 0.
\end{align}
\end{subequations}

\noindent The transport equations for $\rho, \theta_{\text{vd}}, m_{\text{v}}, m_{\text{c}}$ all use an SSPRK3 time discretisation. \par
The SIQN method is typically used in GungHo and Gusto with two inner and two outer loops, called the $2\times2$ configuration. For testing the conservative transport scheme, a complication is that the inner loop updates the dry density but not the mixing ratios, so this conserves dry mass but in general will not conserve the tracer mass. Hence, we follow \textcite{bendall2023trilemma} and use SIQN with sixteen outer loops and one inner loop (the $16\times1$ configuration). The sixteen outer iterations enable convergence of the outer loop, meaning that the dry density after the inner loop matches the transported dry density, and hence mass conservation is determined from the transport scheme. As future work, we aim to modify the inner loop to preserve tracer mass, which would enable our scheme to be used with the standard $2\times2$ SIQN configuration. \par
Conservative injection operations are applied to both mixing ratios to represent them in the transport space; this is computed after the projection of the density into the transport space (Table \ref{table:vert_stag_operations}). Like the transport-only test, we transport $\rho, m_{\text{v}},$ and $ m_{\text{c}}$ in $\mathbb{V}_T=\mathbb{\hat{V}}_{\theta}$. We also transport $\theta_{\text{vd}}$ in $\mathbb{V}_T$, although this is not necessary for tracer mass conservation. Physics increments from the saturation adjustment scheme are applied explicitly, using a forward Euler discretisation, after each full SIQN step. Limiters are not required in this test. Further details of the test case setup and SIQN timestepping are given in Appendix \ref{appendix:setup5}. \par
\begin{figure}
    \centering
    \includegraphics[width=0.6\linewidth]{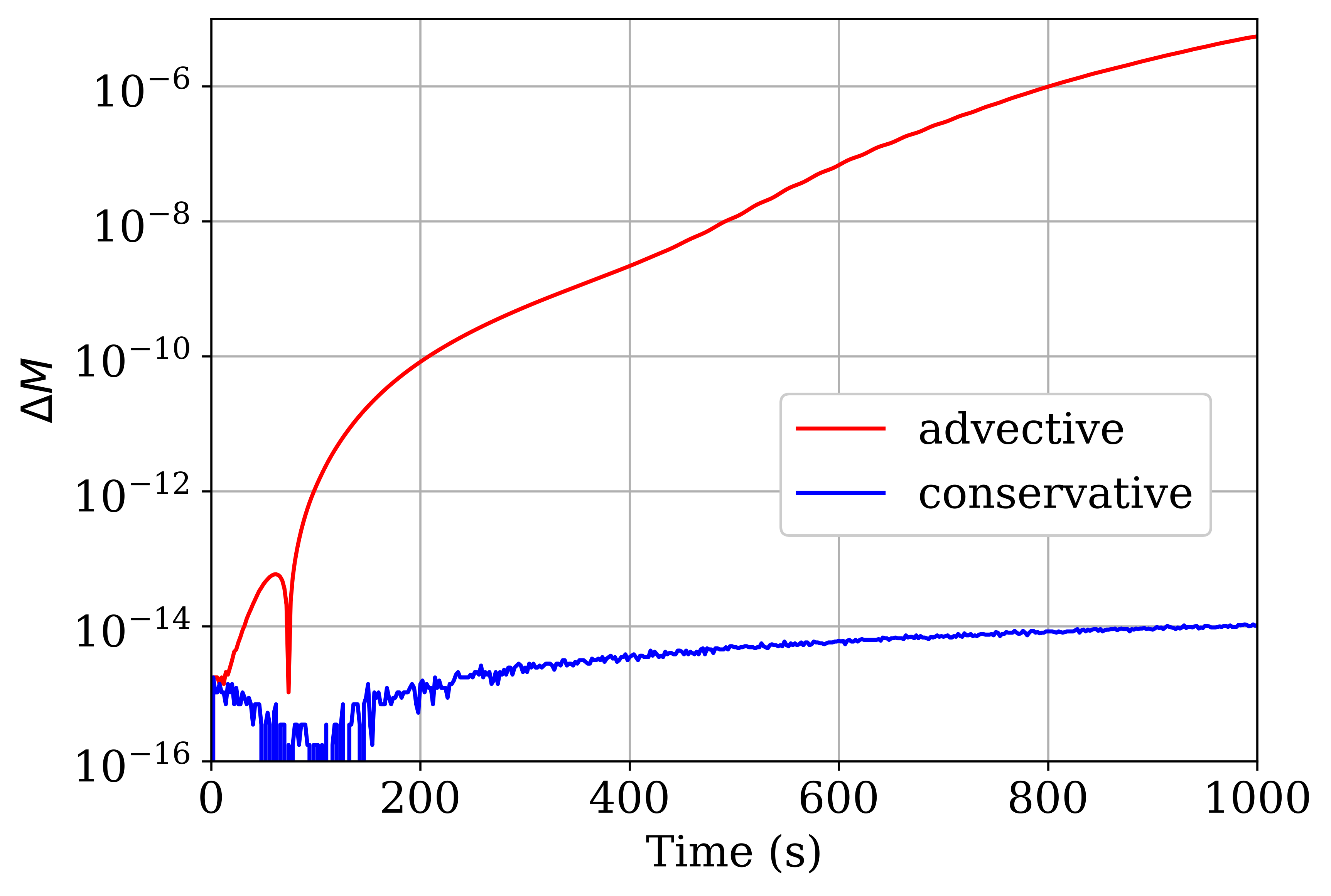}
    \caption{Relative change in total water mass $\Delta M$ (\ref{eq:td_error_over_time}), where $M=\intOmega{\rho(m_{\text{v}}+m_{\text{c}})}$, over time in the Bryan-Fritsch moist rising bubble test. The conservative scheme maintains mass conservation at the level of round-off error.}
    \label{fig:bryan_fritsch_Td}
\end{figure}

Figure \ref{fig:bryan_fritsch_comps} compares the use of advective and conservative tracer transport on the wet equivalent potential temperature $\theta_e$ \parencite{durran1982effects,bryan2002benchmark}. The final $\theta_e$ fields from both schemes are extremely similar. However, Figure \ref{fig:bryan_fritsch_Td} shows there is a marked difference in mass conservation; the advective scheme leads to a growing mass discrepancy over time, whilst the water mass is conserved with the conservative scheme. \par

\begin{figure}
    \centering
    \includegraphics[width=\linewidth]{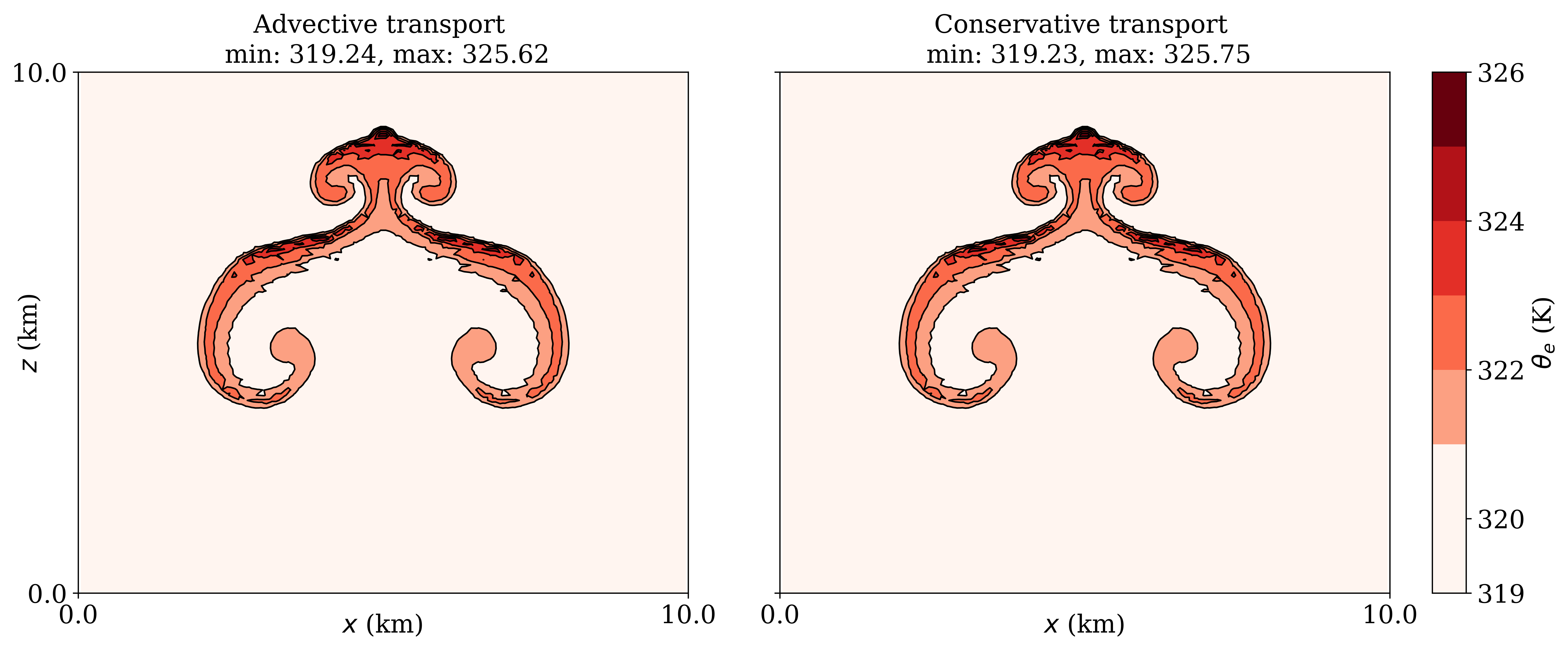}
    \caption{Plots of the wet equivalent potential temperature ($\theta_e$) at the end of the Bryan-Fritsch moist rising bubble simulations with the advective (left) and conservative (right) transport schemes. 1 K contours are used, with the 320 K contour omitted to avoid emphasising small variations around the reference state. The fields are visually indistinguishable between the two transport schemes.}
    \label{fig:bryan_fritsch_comps}
\end{figure}

\clearpage

\section{Conclusions}
\label{section:conclusions}
This paper presented a conservative, discontinuous Galerkin, tracer transport scheme within a compatible finite element discretisation. The novelty of the method is its ability to work with a mixing ratio that is either co-located or vertically staggered from the dry density. In the vertically staggered case, injection and Galerkin projection operations are used to move the density and mixing ratios between their original function space and the shared transport space. To ensure mass conservation when mapping a mixing ratio, we introduce conservative Galerkin projection and conservative injection operators that act upon the tracer density of $\rho m$. These conservative operations use three input fields: the mixing ratio and density in their original function spaces, and the density in the mapped function space. The conservative projection operator was made consistent by subtracting a global mean of the mixing ratio. \par 
Simulations using the Gusto code library were used to validate our conservative transport scheme. Mass conservation and consistency were shown in transport-only tests with co-located and vertically staggered function spaces. The conservative scheme led to the same convergence rate as the advective scheme, whilst allowing for mass conservation near machine precision. Consistency was shown through the preservation of spatially constant mixing ratio fields after transport. Monotonicity, in the co-located case, was addressed through a mean mixing ratio limiter, which blends a next-to-lowest-order mixing ratio with a lowest-order counterpart to ensure mass conservation is maintained. The conservative transport scheme and mean mixing ratio limiter preserve linear correlations, which was seen in the terminator toy test. A final dynamical core-like test with the Bryan-Fritsch rising bubble showed that the conservative scheme can successfully conserve the moist mass in the compressible Euler equations when using a vertical staggering of the mixing ratios from the dry density. \par
A next step is extending the mean mixing ratio limiter to work with the vertically staggered configuration; this requires limiting fields that have continuity requirements in the vertical dimension. We also aim to extend the global conservation of the vertically staggered scheme to local conservation, which will likely require the shifted mesh of \textcite{bendall2023trilemma}. Developing a mass-conserving inner loop for SIQN is also desirable, as this means that the conservative transport equation can ensure mass conservation in the default $2\times2$ configuration. We would also like to apply the conservative tracer transport scheme to other equation sets and configurations such as the family of moist shallow water equations from \textcite{hartney2025exploring}. 






\appendix

\section{Mean mixing ratio proofs}
This appendix provides proofs about the mean mixing ratio (MMR).

\subsection{Monotonicity of the MMR}
\label{append:monotonic_mmr}
Here, we prove monotonicity of the MMR, provided a CFL-type constraint is obeyed, which ensures that the MMR limiter of Section \ref{subsec:monotonicity} enforces monotonicity. The MMR applies to the co-located case of $\rho,m\in\mathbb{V}_{\rho}$. The proofs will take $\mathbb{V}_{\rho}$ to have its degrees of freedom at the vertices of each cell, but the argument holds for any representation of this fully discontinuous function space. We first introduce some useful lemmas for the monotonicity proof.

\begin{lemma}
    The cell-integrated values of $\rho$ and $\rho m$ in element $e$ change only due to fluxes through the surface $\partial e$.
\end{lemma}

\begin{pf}
    Consider element $e$ and indicator test functions of
	\begin{equation}
		\gamma, \eta = \left\lbrace
		\begin{matrix}
			1, & \bm{x}\in e, \\
			0, & \mathrm{otherwise}.
		\end{matrix}\right.
        \label{eq:indicator}
	\end{equation}
	Then, a forward Euler discretisation of the conservative transport system (\ref{eq:mixed_dg_upwind_conserve}) reduces to
	\begin{subequations}
		\begin{align}
			\int_e \widetilde{\rho} \dx{V} & = \int_e \rho \dx{V} - \Delta t \int_{\partial e_{\mathrm{out}}} \rho |\bm{u\cdot}\widehat{\bm{n}}| \dx{S}
			+ \Delta t \int_{\partial e_{\mathrm{in}}} \rho |\bm{u\cdot}\widehat{\bm{n}}| \dx{S}, \label{eq:fv_rho}\\
			\int_e \widetilde{\rho}\widetilde{m} \dx{V} & = \int_e \rho m\dx{V} - \Delta t \int_{\partial e_{\mathrm{out}}} \rho m|\bm{u\cdot}\widehat{\bm{n}}| \dx{S}
			+ \Delta t \int_{\partial e_{\mathrm{in}}} \rho m |\bm{u\cdot}\widehat{\bm{n}}| \dx{S}, \label{eq:fv_tracer}
		\end{align}
    \label{eq:fv_system}
	\end{subequations}
    \noindent with $\widetilde{\rho}, \widetilde{m}$ the updated solutions. The volume integral terms, which involve a gradient of a test function, have vanished. The surface integral terms identify the flow of mass into or out of the element. (\ref{eq:fv_system}) is in the form of a finite volume scheme, and shows that the change in mass (whether of the tracer species or of the reference species) in element $e$ is given by the mass flow into and out of the element. 
\end{pf}

\begin{theorem}[Monotonicity of the forward Euler update MMR]
\label{theorem:monontonic_mmr}
    The MMR in element $e$, $\overline{\widetilde{m}}_e$, computed from the forward Euler update $\widetilde{m}$, is bounded by the prior extrema of the mixing ratio, $m^{i-1,\min}_{e\cup\partial e}$ and $m^{i-1,\max}_{e\cup\partial e}$ (\ref{eq:local_monotonic_bounds}), subject to a CFL-type constraint on the outflow flux, derived below.
\end{theorem}

\begin{pf}
    \renewcommand{\qedsymbol}{}
    Decompose \eqref{eq:fv_rho} and \eqref{eq:fv_tracer} into outflow and inflow parts,
	\begin{equation}
		\int_e \widetilde{\rho} \dx{V} = Q_{\mathrm{out}}+Q_{\mathrm{in}}, \qquad \int_e \widetilde{\rho} \widetilde{m} \dx{V} = R_{\mathrm{out}}+R_{\mathrm{in}}.
	\end{equation}
    \noindent The outflow quantities only draw on values internal to element $e$,
	\begin{equation}
    Q_{\mathrm{out}} := \int_e \rho \dx{V} - \Delta t \int_{\partial e_{\mathrm{out}}} \rho |\bm{u\cdot}\widehat{\bm{n}}| \dx{S}, \quad
		R_{\mathrm{out}} := \int_e \rho m \dx{V} - \Delta t \int_{\partial e_{\mathrm{out}}} \rho m |\bm{u\cdot}\widehat{\bm{n}}| \dx{S} \label{eqn:Tout_Mout}.
	\end{equation}
    \noindent The inflow quantities only draw on values from neighbouring cells,
	\begin{equation}
        Q_{\mathrm{in}} := \Delta t \int_{\partial e_{\mathrm{in}}} \rho |\bm{u\cdot}\widehat{\bm{n}}| \dx{S}, \qquad
		R_{\mathrm{in}} := \Delta t \int_{\partial e_{\mathrm{in}}} \rho m |\bm{u\cdot}\widehat{\bm{n}}| \dx{S}, \label{eqn:Tin_Min}
	\end{equation}
	\noindent with $\rho,m$ the upwind values on the facet $\partial e_{\mathrm{in}}$. \par
    To complete the proof, we provide bounds on $R_{\mathrm{in}}$ and $R_{\mathrm{out}}$ through the following two lemmas. We note that the outflow lemma (Lemma \ref{lemma:outflow}) is similar to Theorem 1 of \textcite{guba2014optimization}. 
\end{pf}

\begin{lemma}[Inflow]
\label{lemma:inflow}
    \textit{Provided $\rho \geq 0$ on the neighbouring side of $\partial e_{\mathrm{in}}$, we have the following bounds on the inflow terms:}
	\begin{equation}
    		Q_{\mathrm{in}}\geq0, \qquad m^{i-1,\min}_{\partial e_{\mathrm{in}}} Q_{\mathrm{in}} \leq R_{\mathrm{in}} \leq m^{i-1,\max}_{\partial e_{\mathrm{in}}} Q_{\mathrm{in}}.
	\end{equation}
\end{lemma}

\begin{pf}
    On $\partial e_{\mathrm{in}}$, we have the bound $m\in \left[m^{i-1,\min}_{\partial e_{\mathrm{in}}},m^{i-1,\max}_{\partial e_{\mathrm{in}}} \right]$. Given that $\rho \geq 0$, we have that $\rho |\bm{u\cdot}\widehat{\bm{n}}| \geq0$ and $Q_{\mathrm{in}}\geq0$ by \eqref{eqn:Tin_Min}. Multiplying the bound on $m$ by $\rho|\bm{u\cdot}\widehat{\bm{n}}|$ gives
	\begin{equation}
		m^{i-1,\min}_{\partial e_{\mathrm{in}}}\,\rho |\bm{u\cdot}\widehat{\bm{n}}| \ \leq\ \rho m |\bm{u\cdot}\widehat{\bm{n}}| \ \leq\ m^{i-1,\max}_{\partial e_{\mathrm{in}}}\,\rho |\bm{u\cdot}\widehat{\bm{n}}| \qquad \text{on } \partial e_{\mathrm{in}},
	\end{equation}
	and integrating over $\partial e_{\mathrm{in}}$ preserves the inequality, giving the result for $R_{\mathrm{in}}$ via \eqref{eqn:Tin_Min}. There is no restriction relating to $\Delta t$ on either inflow term. 
\end{pf}

\begin{lemma}[Outflow]
\label{lemma:outflow}
    Provided $\rho \geq 0$ and that at every vertex $j$ of $e$ on $\partial e_{\mathrm{out}}$ (with $\phi_j$ the nodal basis function of vertex $j$), we have the CFL-type condition
    \begin{equation}
		V^e_j - \Delta t F^{\mathrm{out}}_j \geq 0, \label{eqn:cfl}
	\end{equation}
    \noindent where
    \begin{equation}
       V^e_j := \int_e \rho \phi_j \dx{V}, \quad \quad F^{\mathrm{out}}_j := \int_{\partial e_{\mathrm{out}}} \rho \phi_j |\bm{u\cdot}\widehat{\bm{n}}| \dx{S}. \label{eq:ve_fout_defs}
    \end{equation}
    \noindent We then have the following bounds on the outflow terms:
    \begin{equation}
		Q_{\mathrm{out}} \geq 0, \qquad m^{i-1,\min}_e Q_{\mathrm{out}} \leq R_{\mathrm{out}} \leq m^{i-1,\max}_e Q_{\mathrm{out}}. 
	\end{equation}
\end{lemma}

\begin{pf}
    In element $e$, write $m = (1-\nu)m^{i-1,\min}_e + \nu m^{i-1,\max}_e$, with $\nu\in\mathbb{V}_T$ and the constraint that $\nu\in[0,1]$. Using \eqref{eqn:Tout_Mout}, the outflow terms become
	\begin{equation}
		Q_{\mathrm{out}}=X(1)=X(\nu)+X(1-\nu), \qquad R_{\mathrm{out}} = m^{i-1,\min}_e X(1-\nu) + m^{i-1,\max}_e X(\nu), \label{eq:Q_M_outflow}
	\end{equation}
	where
	\begin{equation}
	X(\nu) := \int_e \rho \nu \dx{V} - \Delta t\int_{\partial e_{\mathrm{out}}} \rho \nu |\bm{u\cdot}\widehat{\bm{n}}| \dx{S}.
	\end{equation}
	Using a nodal basis for $\mathbb{V}_T$, $\nu=\sum_j \nu_j\phi_j$ with the constraint $\nu\in[0,1]$ throughout $e$ holding if and only if $\nu_j\in[0,1]$ at each vertex $j$, since $\phi_j\geq0$ and $\sum_j\phi_j=1$; the same applies to $1-\nu$, whose nodal values are $1-\nu_j\in[0,1]$.
	\\
	\\
	Since $X$ is linear,
	\begin{equation}
		X(\nu) = \sum_j \nu_j X(\phi_j), \qquad X(1-\nu) = \sum_j (1-\nu_j) X(\phi_j), \label{eq:X_sums}
	\end{equation}
	where $X(\phi_j) = V^e_j - \Delta t F^{\mathrm{out}}_j$ directly from the definitions of $V^e_j$ and $F^{\mathrm{out}}_j$ (\ref{eq:ve_fout_defs}). So, requiring $X(\nu)\geq0$ for every admissible $\nu$ is equivalent to requiring $X(\phi_j)\geq0$ at every vertex $j$, which is exactly the CFL-type constraint $V^e_j - \Delta t F^{\mathrm{out}}_j \geq 0$ \eqref{eqn:cfl}. \par
    Given \eqref{eqn:cfl}, every term $\nu_j X(\phi_j)$ and $(1-\nu_j)X(\phi_j)$ in the two sums above is a product of two non-negative factors, so
	\begin{equation}
		X(\nu)\geq0, \qquad X(1-\nu)\geq0.
	\end{equation}
	In particular, $Q_{\mathrm{out}}=X(\nu)+X(1-\nu)\geq0$, and $R_{\mathrm{out}}$ is a weighted average of $m^{i-1,\min}_e$ and $m^{i-1,\max}_e$ with these same non-negative weights: since $m^{i-1,\max}_e\geq m^{i-1,\min}_e$,
	\begin{equation}
		m^{i-1,\min}_e Q_{\mathrm{out}} = m^{i-1,\min}_e\left[X(1-\nu)+X(\nu)\right] \leq m^{i-1,\min}_e X(1-\nu) + m^{i-1,\max}_e X(\nu) = R_{\mathrm{out}},
	\end{equation}
	\begin{equation}
		R_{\mathrm{out}} = m^{i-1,\min}_e X(1-\nu) + m^{i-1,\max}_e X(\nu) \leq m^{i-1,\max}_e\left[X(1-\nu)+X(\nu)\right] = m^{i-1,\max}_e Q_{\mathrm{out}}. 
	\end{equation}
\end{pf}

\begin{pf}[Monotonicity of the mean mixing ratio]
	Combining the bounds on $R_{\text{out}}$ and $R_{\text{in}}$ gives
	\begin{equation}
		m^{i-1,\min}_{e\cup\partial e}\left(Q_{\textrm{out}}+Q_{\textrm{in}}\right) \leq R_{\text{out}}+R_{\text{in}} \leq m^{i-1,\max}_{e\cup\partial e}\left(Q_{\textrm{out}}+Q_{\textrm{in}}\right).
	\end{equation}
	Since it is a requirement that $Q_{\textrm{out}}+Q_{\textrm{in}} =\int_e\widetilde{\rho} \dx{V}> 0$ (\ref{eq:positive_rho}), we can safely divide the bound through by $Q_{\textrm{out}} + Q_{\textrm{in}}$ to give the desired bound
	\begin{equation}
		m^{i-1,\min}_{e\cup\partial e} \leq \overline{\widetilde{m}}_e \leq m^{i-1,\max}_{e\cup\partial e}. 
	\end{equation}
\end{pf}

We have shown that the MMR is monotonic after a forward Euler step of the transport equation. However, the limiter applies to the stage $i$ field, not the forward Euler update. It is a simple extension to show that the stage $i$ MMR is also monotonic.

\begin{corollary}[Monotonicity of the stage $i$ MMR]
\label{corollary:stage_mmr_monotonic}
    Given a monotonic MMR after each forward Euler step, $\overline{\widetilde{m}}$, the stage $i$ MMR is also monotonic and obeys
    \begin{equation}
		m^{i,-}_{e} \leq \overline{m}^{i}_e \leq m^{i,+}_{e},
	\end{equation}
    \noindent where $m^{i,-}_{e},m^{i,+}_{e}$ are the local bounds given by (\ref{eq:high_low_bounds}).
\end{corollary}

\begin{pf}
    The mean mixing ratio (\ref{eq:mmr}) in element $e$ is given by
    \begin{equation}
        \overline{m}^{i}_e
        = \frac{\int_e\rho^{i}m^{i}\dx{V}}{\int_e\rho^{i}\dx{V}}.
    \end{equation}
    \noindent Now, consider a general SSPRK3 stage (\ref{eq:ssprk3_gen_stage}) and integrate it over element $e$; this is valid as the stage is solved in weak form and the indicator function (\ref{eq:indicator}) lies in the transport space $\mathbb{V}_T$:
    \begin{equation}
    \int_e\rho^{i}\dx{V} = \alpha_i\int_e\rho^{n}\dx{V} + (1-\alpha_i)\int_e\widetilde\rho\dx{V}, \qquad
    \int_e\rho^{i}m^{i}\dx{V} = \alpha_i\int_e\rho^{n}m^{n}\dx{V} + (1-\alpha_i)\int_e\widetilde\rho\widetilde m\dx{V},
\end{equation}
    This means that
    \begin{equation}
        \overline{m}^{i}_e
        = \frac{\int_e\rho^{i}m^{i}\dx{V}}{\int_e\rho^{i}\dx{V}} 
        = \alpha_i\frac{\int_e\rho^{n}m^{n}\dx{V}}{\int_e\rho^{i}\dx{V}} + (1-\alpha_i) \frac{\int_e\widetilde\rho\widetilde m\dx{V}}{\int_e\rho^{i}\dx{V}}
        = w\,\overline{m}^n_e + (1-w)\,\overline{\widetilde m}_e,
    \label{eq:mmr_convex_combo}
    \end{equation}
    with weights
    \begin{equation}
        w = \frac{\alpha_i\int_e\rho^n\dx{V}}{\int_e\rho^{i}\dx{V}} = \frac{\alpha_i\int_e\rho^n\dx{V}}{\alpha_i\int_e\rho^n\dx{V}+(1-\alpha_i)\int_e\widetilde\rho\dx{V}},
    \end{equation}
    and
    \begin{equation}
        1-w = \frac{\int_e \rho^i \dx{V} -\alpha_i\int_e\rho^n\dx{V}}{\int_e\rho^{i}\dx{V}} = \frac{(1-\alpha_i)\int_e\widetilde\rho\dx{V}}{\int_e\rho^{i}\dx{V}}.
    \end{equation}

    As $\int_e\widetilde\rho\dx{V}>0$ and $\int_e\rho^n\dx{V}>0$ (\ref{eq:positive_rho}), and $\alpha_i \in [0,1)$, it follows that $w\in[0,1]$, and so $\overline{m}^{i}_e$ (\ref{eq:mmr_convex_combo}) is the convex combination of $\overline{m}^n_e$ and $\overline{\widetilde m}_e$. As $m^{n,\min}_e \leq \overline{m}^n_e \leq m^{n,\max}_e$ and $m^{i-1,\min}_{e\cup\partial e} \leq \overline{\widetilde m}_e \leq m^{i-1,\max}_{e\cup\partial e}$ by Theorem \ref{theorem:monontonic_mmr}, we have the common bound $\overline{m}^n_e, \overline{\widetilde m}_e \in [m^{i,-}_{e},m^{i,+}_{e}]$, so it follows that $\overline{m}^{i}_e \in [m^{i,-}_{e},m^{i,+}_{e}]$. \par
    
\end{pf}

\subsection{The maximal blending function}
\label{append:proof_blending_func}

We now show that using the maximal blending function (\ref{eq:multiple_blend_func}) to limit multiple mixing ratios will preserve linear correlations. We consider the case of the linear correlation forming a spatially constant field, as per the terminator toy test (Section \ref{subsec:terminator}).

\begin{theorem}
\label{theorem:maximal_blending_fun}
    The maximal blending function $\Lambda_e^M$, used to limit all mixing ratios $m_k$ in element $e$, preserves a linear correlation of
    \begin{equation}
        \sum_{k=1}^N a_k m_k = D,
    \end{equation}
    \noindent where $a_k \in \mathbb{R}$ and $D$ is a conserved, spatially constant field.
\end{theorem}

\begin{pf}
    Consider the linear combination in element $e$ after transport and limiting, with $\Lambda_{k,e}$ the blending weight (\ref{eq:lambda_e}) for tracer $k$,
    \begin{align}
        \sum_{k=1}^N a_k m_{k,e}^* &= \sum_{k=1}^N a_k [(1-\Lambda_{k,e})m_{k,e} + \Lambda_{k,e} \overline{m}_{k,e}] , \\
        &= \sum_{k=1}^N a_k m_{k,e} - \sum_{k=1}^N a_k \Lambda_{k,e} m_{k,e} + \sum_{k=1}^N a_k \Lambda_{k,e} \frac{\int_e \rho m_{k,e} \dx{V}}{\int_e \rho \dx{V}}, \\
        &= D + \left[ \frac{\int_e \rho \sum_{k=1}^N \Lambda_{k,e} a_k m_{k,e} \dx{V}}{\int_e \rho \dx{V}} - \sum_{k=1}^N a_k \Lambda_{k,e} m_{k,e} \right].
    \end{align}
    \noindent Note that $\sum_{k=1}^N a_k m_{k,e}=D$ follows from linearity and consistency of the transport update (Theorem \ref{theorem:consistency}). 
    If we use the maximal blending function for all mixing ratios, then $\sum_{k=1}^N a_k \Lambda_{k,e} m_{k,e} = \Lambda_e^M \sum_{k=1}^N a_k m_{k,e} = \Lambda_e^M D$, and so the bracketed term simplifies to zero, leaving $\sum_{k=1}^N a_k m_{k,e}^*=D$. 
\end{pf}

\section{Test case setups}
This appendix provides details about the test cases used in this paper. The description of the co-located and vertical slice domains and their boundary conditions can be found at the beginning of Section \ref{subsec:scheme_props}.

\subsection{Co-located transport-only test}
\label{appendix:setup1}
We use the divergent flow field described in \textcite{lauritzen2012standard}. Using $u$ to denote the zonal flow along lines of constant latitude, $\phi \in [-\pi/2, \pi/2]$, and $v$ the meridional flow along lines of constant longitude, $\lambda \in [-\pi, \pi)$, the velocity field is
\begin{subequations}
\begin{align}
    u(\lambda, \phi, t) &= - \frac{5a}{\tau} \sin^2 \bracfrac{\lambda'}{2}  \sin(2\phi) \cos^2(\phi) \cos \bracfrac{\pi t}{\tau} + \frac{2 \pi a}{\tau} \cos(\phi), \label{eq:colocated_flow_field_u} \\ 
    v(\lambda, \phi, t) &= \frac{5a}{2\tau} \sin(\lambda') \cos^3(\phi) \cos \bracfrac{\pi t}{\tau},
\end{align}
\label{eq:colocated_flow_field}
\end{subequations}

\noindent where $a=6371.22$ km is taken as the Earth's radius, $\lambda^{'} = \lambda - 2 \pi t/\tau$, and $\tau$ is the time period for the flow to return to its initial state, which is set to the simulation length of twelve days ($\tau = 1,036,800$ s). The mean flow term of $2 \pi a \cos(\phi) / \tau$ in (\ref{eq:colocated_flow_field_u}) breaks the symmetry of the flow field to ensure that errors generated in the first half of the simulation are not reversed in the second half, which can lead to misleadingly small transport errors. A timestep size of $\Delta t = 450$ s is used for all tests. \par
The smooth scalar field option from \textcite{nair2010class} is chosen to initialise the mixing ratio in the convergence test and the density in the consistency test. This defines Gaussian perturbations of
\begin{equation}
    g_i(\lambda_i, \phi_i) = g_{\text{max}} \exp (-b_0[(X - X_i)^2 + (Y - Y_i)^2 + (Z - Z_i)^2] ),
\end{equation}

\noindent where $(X,Y,Z)$ identify the Cartesian coordinate on the unit sphere, computed from longitude-latitude $(\lambda,\phi)$ coordinates using
\begin{equation}
    X = \cos(\phi) \cos(\lambda), ~Y = \cos(\phi) \sin(\lambda), ~Z = \sin(\phi).
\end{equation}

\noindent We select a Gaussian decay rate of $b_0 = 5$ and centre the two Gaussian bumps at $(\lambda_1,\phi_1)=(-\pi/4,0)$ and $(\lambda_2,\phi_2)=(\pi/4,0)$. \par
The convergence setting uses
\begin{subequations}
    \begin{align}
        \rho(\lambda,\phi,t=0) &= \rho_b + 0.5 \cos(\phi), \label{eq:colocated_test_conv_rho} \\
        m(\lambda,\phi,t=0) &= m_0 + g_1 + g_2,
    \end{align}
\end{subequations}

\noindent with $\rho_b = 1 ~\text{kg} ~\text{m}^{-3},  m_0 = 0.02 ~\text{kg} ~\text{kg}^{-1}$, and $g_{\text{max}} = 0.05 ~\text{kg} ~\text{kg}^{-1}$. Spatial resolutions of $N_e \in \{ 16,24,32,48\} $ are used, where $N_e$ denotes the number of cells along a cubed-sphere panel edge. The largest Courant number tested is approximately 0.2. \par
The consistency setting uses
\begin{subequations}
    \begin{align}
        \rho(\lambda,\phi,t=0) &= \rho_b + g_1 + g_2, \label{eq:colocate_consist_rho} \\ 
        m(\lambda,\phi,t=0) &= m_0,
    \end{align}
\end{subequations}

\noindent with $\rho_b = 1 ~\text{kg} ~\text{m}^{-3}, ~m_0 = 0.02 ~\text{kg} ~\text{kg}^{-1}$, $g_{\text{max}} = 0.5 ~\text{kg} ~\text{m}^{-3}$. The spatial resolution is $N_e = 32$. \par

\subsection{Co-located MMR limiter test}
\label{appendix:setup2}
The test for the mean mixing ratio limiter is similar to the convergence configuration of the co-located transport-only test, and uses the same prescribed velocity (\ref{eq:colocated_flow_field}) and initial density field (\ref{eq:colocated_test_conv_rho}). The mixing ratio is initialised by two slotted-cylinder functions \parencite{nair2010class},
\begin{equation}
    m(\lambda, \phi, t=0) = 
    \begin{cases}
    1, \text{if} ~r_1 < 0.5, |\lambda - \lambda_1| > 1/12, \\
    1, \text{if} ~r_2 < 0.5, |\lambda - \lambda_2| > 1/12,  \\
    1, \text{if} ~r_1 < 0.5, |\lambda - \lambda_1| \leq 1/12, (\phi-\phi_1) < -5/24, \\
    1, \text{if} ~r_2 < 0.5, |\lambda - \lambda_2| \leq 1/12, (\phi-\phi_2) > 5/24, \\
    0, ~\text{else},
    \end{cases}
\end{equation}

\noindent with the cylinders centred at $(\lambda_1,\phi_1) = (-\pi/4, 0)$ and $(\lambda_2,\phi_2) = (\pi/4, 0)$, and $r_i$ the great arc angle (great circle distance on the unit sphere) of
\begin{equation}
    r_i(\lambda,\phi) = \arccos ( \sin(\phi) \sin(\phi_i) + \cos(\phi) \cos(\phi_i) \cos(\lambda - \lambda_i)).
\end{equation}

A spatial resolution of $N_e = 48$ cells per panel edge and a timestep of $\Delta t = 450$ s are used. The maximum Courant number in this test is approximately 0.2.

\subsection{Co-located terminator toy test}
\label{appendix:setup3}
The terminator toy test is described in \textcite{lauritzen2015terminator}. The reaction rates for the interaction of chemical species, $k_1$ and $k_2$, are given by

\begin{subequations}
\begin{align}
    k_1(\lambda,\phi) &= \max \{ 0, \sin(\phi) \sin(\phi_{\text{c}}) + \cos(\phi) \cos(\phi_{\text{c}}) \cos(\lambda - \lambda_{\text{c}}) \}, \\
    k_2(\lambda,\phi) &= 1,
\end{align}
\label{eq:k_rates_term_toy}
\end{subequations}

\noindent with $(\lambda_{\text{c}}, \phi_{\text{c}}) = ( -\pi / 3, \pi/9)$ the location of the peak reaction rate of $k_1 = 1$. \par
The initial, spatially constant, total mixing ratio of $X_T (0) = 4 \times 10^{-6} ~\text{kg} ~\text{kg}^{-1}$ is distributed between the two species by
\begin{equation}
    X(0) = d - r, ~ X_2(0) = 0.5(X_T(0) - d + r),
\end{equation}

\noindent where $r=k_1 / (4 k_2),$ and $ d = \sqrt{r^2 + 2r X_T}$. The initial condition for the dry density is the Gaussian bumps used in the co-located consistency test (\ref{eq:colocate_consist_rho}). The prescribed non-divergent velocity field is
\begin{subequations}
\begin{align}
    u(\lambda, \phi, t) &=  \frac{10a}{\tau} \sin^2 (\lambda^{'}) \sin(2 \phi) \cos \bracfrac{\pi t}{\tau} + \frac{2 \pi a}{\tau} \cos(\phi), \\
    v(\lambda, \phi, t) &= \frac{10a}{\tau} \sin(2 \lambda') \cos(\phi) \cos \bracfrac{\pi t}{\tau},
\end{align}
\end{subequations}

\noindent with $\lambda^{'} = \lambda - 2 \pi t/\tau$. The full simulation time is twelve days, i.e. $\tau = 1,036,800$ s. A spatial resolution of $N_e = 24$ cells per panel edge and a timestep of $\Delta t = 450$ s are used. The maximum Courant number in this test is approximately 0.2. \par

\subsection{Vertically staggered transport-only test}
\label{appendix:setup4}
The transporting velocity of $\vec{u} = [u, w]^{\text{T}}$, with $u$ the horizontal component of the wind and $w$ the vertical component of the wind, is described by
\begin{subequations}
\begin{align}
    u(x,z,t) &= U - \frac{W \pi L_x}{H_z} \cos \left( \frac{\pi t}{\tau} \right) \cos \left( \frac{2 \pi x'}{L_x} \right) \cos \left( \frac{\pi z }{H_z} \right), \\
    w(x,z,t) &= 2 \pi W \cos \left( \frac{\pi t}{\tau} \right) \sin \left( \frac{2 \pi x'}{L_x} \right) \sin \left( \frac{\pi z}{H_z} \right),
\end{align}
\end{subequations}

\noindent with $U = L_x/\tau, W = U/10, x' = x - Ut$. The domain is a square of length $L_x = 2000 ~\text{m}$ and height $ H_z = 2000 ~\text{m}$. The cells are set to a unity aspect ratio of $\Delta x = \Delta z = L_x/N$, where $N$ is the number of elements in each dimension. The simulation time is $\tau = 2000$ s, and a timestep of $\Delta t=2 ~\text{s}$ is used for all tests.  \par
Gaussian perturbations are defined on the vertical slice domain as
\begin{equation}
    f(x_c,z_c) = f_0 \exp( -l^2(x_c, z_c) / l_c^2),
\end{equation}

\noindent with $f_0$ the amplitude of the perturbation,  $l^2(x_c,z_c) = (\min\{ |x-x_c|, L_x - |x-x_c| \})^2 + |z-z_c|^2$ the shortest squared distance from the centre of a perturbation at ($x_c, z_c$), and $l_c = 2 L_x/25$ the decay rate of the Gaussian functions. The Gaussian perturbations are centred at $(x_1, z_1) = (3L_x/8, H_z/2)$ and $(x_2, z_2)=(5L_x/8, H_z/2)$. \par
The convergence setting uses
\begin{subequations}
\begin{align}
    \rho(x,z,t=0) &= \rho_b + (\rho_t - \rho_b) \frac{z}{H_z}, \\
    m(x,z,t=0) &= m_0 + f(x_1, z_1) + f(x_2, z_2),
\end{align}
\end{subequations}

\noindent with $\rho_b = 1 ~\text{kg} ~\text{m}^{-3}$, $\rho_t = 0.5 ~\text{kg} ~\text{m}^{-3}$, $m_0 = 0.02 ~\text{kg} ~\text{kg}^{-1}$, and $f_0 = 0.05 ~\text{kg} ~\text{kg}^{-1}$. The spatial resolutions are $N \in \{ 50, 60, 70, 80, 90, 100 \}$, with a largest tested Courant number of approximately 0.1. \par

The consistency setting uses
\begin{subequations}
\begin{align}
    \rho(x,z,t=0) &= \rho_b + f(x_1, z_1) + f(x_2, z_2), \\
    m(x,z,t=0) &= m_0,
\end{align}
\end{subequations}

\noindent with $\rho_b = 0.5 ~\text{kg} ~\text{m}^{-3}$, $f_0 = 0.5 ~\text{kg} ~\text{m}^{-3}$, and $m_0 = 0.02 ~\text{kg} ~\text{kg}^{-1}$. The spatial resolution is $N=100$. \par

\subsection{Vertically staggered dynamical core test}
\label{appendix:setup5}
The test parameters and initial condition follow those in \textcite{bendall2020compatible}, which is a vertical slice modification of the rising bubble test from \textcite{bryan2002benchmark}. The domain is of length $L_x = 10$ km and height $H_z = 10$ km. Fifty elements are used along both the length and the height dimensions, so $\Delta x = \Delta z = 200$ m. The simulation length is 1000 s and the timestep size is $\Delta t =2$ s. The maximum Courant number is 0.175. \par
As noted in \textcite{bendall2023trilemma}, the constant mixing ratio initialisation of \textcite{bryan2002benchmark} can lead to surprisingly good conservation with advective form transport; this is due to a cancellation of individual conservation errors in the cloud water and water vapour fields in the tracer density of $\rho_X=\rho(m_{\text{v}} + m_{\text{c}})$. This fortuitous mass conservation with advective transport was indeed observed in initial tests. To avoid the cancellation of mixing ratio errors, the initial total mixing ratio ($m=m_{\text{v}} + m_{\text{c}}$) is defined to vary linearly in the vertical coordinate,
\begin{equation}
    m(x,z,t=0) = m_0 - \frac{z \delta m}{H_z},
\end{equation}

\noindent with $m_0 = 0.02 ~\text{kg} ~\text{kg}^{-1}, \delta m = 0.005 ~\text{kg} ~\text{kg}^{-1}$. The initial water vapour is set to a saturation value, with the cloud water comprising the remainder of the initial mixing ratio. A perturbation is applied to the initial $\theta_{\text{vd}}$ field. \par
A novel feature of the SIQN method in Gusto, compared to GungHo, is the ability to use different time discretisations for the transport of each field in the outer loop. We use the implicit midpoint rule for the transport of $\vec{u}$ and SSPRK3 for the other fields ($\rho,\theta_{\text{vd}},m_{\text{v}},m_{\text{c}}$). In each iteration of the SIQN outer loop, an advecting velocity is first identified, then used to transport each field. The transport of each field is typically performed separately in Gusto, but this was modified to allow $\rho,m_{\text{v}},$ and $m_{\text{c}}$ to be transported simultaneously in the conservative tracer scheme. \par

\section*{Code Availability}
\sloppy
The GitHub repository for Gusto can be found at \mbox{\url{https://github.com/firedrakeproject/gusto}}. Python scripts for running the Gusto simulations and figure making can be found at the GitHub repository \mbox{\url{https://github.com/ta440/gusto_conservative_transport}}. Many of the visualisations used plotting functions from the tomplot GitHub repository \mbox{\url{https://github.com/tommbendall/tomplot}}.

\section*{Declaration of generative AI and AI-assisted technologies in the manuscript preparation process}
During the preparation of this work, the authors used Claude Opus 5.5 to implement the monotonic aspects of the mean mixing ratio limiter, give feedback on the proofs, and check the manuscript for typographical errors. After using this tool, the authors reviewed and edited the content as needed and take full responsibility for the content of the published article.

\printbibliography

@article{cotter2016embedded,
  title={Embedded discontinuous {G}alerkin transport schemes with localised limiters},
  author={Cotter, Colin J and Kuzmin, Dmitri},
  journal={J. Comput. Phys.},
  volume={311},
  pages={363--373},
  year={2016},
  publisher={Elsevier},
  doi={10.1016/j.jcp.2016.02.021}
}

@article{bendall2019recovered,
  title={The `recovered space’ advection scheme for lowest-order compatible finite element methods},
  author={Bendall, Thomas M and Cotter, Colin J and Shipton, Jemma},
  journal={J. Comput. Phys.},
  volume={390},
  pages={342--358},
  year={2019},
  publisher={Elsevier},
  doi={10.1016/j.jcp.2019.04.013}
}

@article{bendall2023trilemma,
  title={A solution to the trilemma of the moist {C}harney--{P}hillips staggering},
  author={Bendall, Thomas M and Wood, Nigel and Thuburn, John and Cotter, Colin J},
  journal={Quart.\ J.\ Roy.\ Meteor.\ Soc.},
  volume={149},
  number={750},
  pages={262--276},
  year={2023},
  publisher={Wiley Online Library},
  doi={10.1002/qj.4406}
}

@article{Staniforth_hor_grids,
  title={Horizontal grids for global weather and climate prediction models: a review},
  author={Staniforth, Andrew and Thuburn, John},
  journal={Quart.\ J.\ Roy.\ Meteor.\ Soc.},
  volume={138},
  number={662},
  pages={1--26},
  year={2012},
  publisher={Wiley Online Library},
  doi={10.1002/qj.958}
}

@article{adams2019lfric,
  title={{LFR}ic: Meeting the challenges of scalability and performance portability in Weather and Climate models},
  author={Adams, Samantha V and Ford, Rupert W and Hambley, M and Hobson, JM and Kav{\v{c}}i{\v{c}}, I and Maynard, Christopher M and Melvin, Thomas and M{\"u}ller, Eike Hermann and Mullerworth, S and Porter, Andrew R and others},
  journal={J. Parallel Distrib. Comput.},
  volume={132},
  pages={383--396},
  year={2019},
  publisher={Elsevier},
  doi={10.1016/j.jpdc.2019.02.007}
}

@article{cockburn2001runge,
  title={Runge--{K}utta discontinuous {G}alerkin methods for convection-dominated problems},
  author={Cockburn, Bernardo and Shu, Chi-Wang},
  journal={J. Sci. Comput.},
  volume={16},
  pages={173--261},
  year={2001},
  publisher={Springer},
  doi={10.1023/a:1012873910884}
}

@article{cotter_shipton_mixed_FE,
  title={Mixed finite elements for numerical weather prediction},
  author={Cotter, Colin J and Shipton, Jemma},
  journal={J. Comput. Phys.},
  volume={231},
  number={21},
  pages={7076--7091},
  year={2012},
  publisher={Elsevier},
  doi={10.1016/j.jcp.2012.05.020}
}

@article{melvin2019cartesian,
  title={A mixed finite-element, finite-volume, semi-implicit discretization for atmospheric dynamics: Cartesian geometry},
  author={Melvin, Thomas and Benacchio, Tommaso and Shipway, Ben and Wood, Nigel and Thuburn, John and Cotter, Colin},
  journal={Quart.\ J.\ Roy.\ Meteor.\ Soc.},
  volume={145},
  number={724},
  pages={2835--2853},
  year={2019},
  publisher={Wiley Online Library},
  doi={10.1002/qj.3501}
}

@article{melvin2024mixed,
  title={A mixed finite-element, finite-volume, semi-implicit discretisation for atmospheric dynamics: Spherical geometry},
  author={Melvin, Thomas and Shipway, Ben and Wood, Nigel and Benacchio, Tommaso and Bendall, Thomas and Boutle, Ian and Brown, Alex and Johnson, Christine and Kent, James and Pring, Stephen and others},
  journal={Quart.\ J.\ Roy.\ Meteor.\ Soc.},
  year={2024},
  publisher={Wiley Online Library},
  volume={150},
  number={764},
  pages={4252--4269},
  doi={10.1002/qj.4814}
}

@article{lauritzen2015terminator,
  title={The terminator “toy” chemistry test: A simple tool to assess errors in transport schemes},
  author={Lauritzen, PH and Conley, AJ and Lamarque, J-F and Vitt, F and Taylor, MA},
  journal={Geosci. Model Dev.},
  volume={8},
  number={5},
  pages={1299--1313},
  year={2015},
  publisher={Copernicus GmbH G{\"o}ttingen, Germany},
  doi={10.5194/gmd-8-1299-2015}
}

@article{zerroukat2002slice,
  title={{SLICE}: A {S}emi-{L}agrangian {I}nherently {C}onserving and {E}fficient scheme for transport problems},
  author={Zerroukat, Mohamed and Wood, Nigel and Staniforth, Andrew},
  journal={Quart.\ J.\ Roy.\ Meteor.\ Soc.},
  volume={128},
  number={586},
  pages={2801--2820},
  year={2002},
  publisher={Wiley Online Library},
  doi={10.1256/qj.02.69}
}

@article{bendall2025swift,
  title={{SWIFT}: A Monotonic, Flux-Form Semi-{L}agrangian Tracer Transport Scheme for Flow with Large {C}ourant Numbers},
  author={Bendall, Thomas M and Kent, James},
  journal={Mon.\ Wea.\ Rev.},
  year={2025},
  publisher={American Meteorological Society},
  volume={153},
  number={4},
  pages={565--587},
  doi={10.1175/mwr-d-24-0110.1}
}

@article{thuburn2008some,
  title={Some conservation issues for the dynamical cores of {NWP} and climate models},
  author={Thuburn, John},
  journal={J. Comput. Phys.},
  volume={227},
  number={7},
  pages={3715--3730},
  year={2008},
  publisher={Elsevier},
  doi={10.1016/j.jcp.2006.08.016}
}

@book{Durran,
  title={Numerical methods for fluid dynamics: With applications to geophysics},
  author={Durran, Dale R},
  volume={32},
  year={2010},
  publisher={Springer Science \& Business Media},
  doi={10.1007/978-1-4419-6412-0}
}

@article{nair2010class,
  title={A class of deformational flow test cases for linear transport problems on the sphere},
  author={Nair, Ramachandran D and Lauritzen, Peter H},
  journal={J. Comput. Phys.},
  volume={229},
  number={23},
  pages={8868--8887},
  year={2010},
  publisher={Elsevier},
  doi={10.1016/j.jcp.2010.08.014}
}

@article{cotter2023_acta_numerica,
  title={Compatible finite element methods for geophysical fluid dynamics},
  author={Cotter, Colin J},
  journal={Acta Numerica},
  volume={32},
  pages={291--393},
  year={2023},
  publisher={Cambridge University Press},
  doi={10.1017/s0962492923000028}
}

@article{bryan2002benchmark,
  title={A benchmark simulation for moist nonhydrostatic numerical models},
  author={Bryan, George H and Fritsch, J Michael},
  journal={Mon.\ Wea.\ Rev.},
  volume={130},
  number={12},
  pages={2917--2928},
  year={2002},
  publisher={American Meteorological Society},
  doi={10.1175/1520-0493(2002)130<2917:ABSFMN>2.0.CO;2}
}

@article{herrington2019exploring,
  title={Exploring a lower-resolution physics grid in {CAM}-{SE}-{CSLAM}},
  author={Herrington, Adam R and Lauritzen, Peter H and Reed, Kevin A and Goldhaber, Steve and Eaton, Brian E},
  journal={J. Adv. Model. Earth Syst.},
  volume={11},
  number={7},
  pages={1894--1916},
  year={2019},
  publisher={Wiley Online Library},
  doi={10.1029/2019ms001684}
}

@article{brown2024physics,
  title={Physics--dynamics--chemistry coupling across different meshes in {LFR}ic-{A}tmosphere: {F}ormulation and idealised tests},
  author={Brown, Alex and Bendall, Thomas M and Boutle, Ian and Melvin, Thomas and Shipway, Ben},
  journal={Quart.\ J.\ Roy.\ Meteor.\ Soc.},
  volume={150},
  number={764},
  pages={4650--4670},
  year={2024},
  publisher={Wiley Online Library},
  doi={10.1002/qj.4836}
}

@article{hartney2025exploring,
  title={Exploring forms of the moist shallow-water equations using a new compatible finite-element discretisation},
  author={Hartney, Nell and Bendall, Thomas M and Shipton, Jemma},
  journal={Quart.\ J.\ Roy.\ Meteor.\ Soc.},
  volume={152},
  number={774},
  pages={e70018},
  year={2026},
  publisher={Wiley Online Library},
  doi={10.1002/qj.70018}
}

@article{zhang2008consistency,
  title={Consistency problem with tracer advection in the atmospheric model {GAMIL}},
  author={Zhang, Kai and Wan, Hui and Wang, Bin and Zhang, Meigen},
  journal={Adv. Atmos. Sci.},
  volume={25},
  pages={306--318},
  year={2008},
  publisher={Springer},
  doi={10.1007/s00376-008-0306-z}
}

@article{ronchi1996_equiangular,
  title={The “cubed sphere”: A new method for the solution of partial differential equations in spherical geometry},
  author={Ronchi, Corrado and Iacono, Roberto and Paolucci, Pier S},
  journal={J. Comput. Phys.},
  volume={124},
  number={1},
  pages={93--114},
  year={1996},
  publisher={Elsevier},
  doi={10.1006/jcph.1996.0047}
}

@incollection{balay1997efficient,
  title={Efficient management of parallelism in object-oriented numerical software libraries},
  author={Balay, Satish and Gropp, William D and McInnes, Lois Curfman and Smith, Barry F},
  booktitle={Modern software tools for scientific computing},
  pages={163--202},
  year={1997},
  publisher={Springer},
  doi={10.1007/978-1-4612-1986-6_8}
}

@techreport{balay2024petsc,
  title={{PETS}c/{TAO} Users Manual Revision 3.22},
  author={Balay, Satish and Abhyankar, S and Adams, M and Brown, Jed and Brune, P and Buschelman, K and Constantinescu, EM and Dener, A and Faibussowitsch, J and Gropp, William and others},
  year={2024},
  institution={Argonne National Laboratory (ANL), Argonne, IL (United States)},
  doi={10.2172/2476320}
}

@article{skamarock2012multiscale,
  title={A multiscale nonhydrostatic atmospheric model using centroidal {V}oronoi tesselations and {C}-grid staggering},
  author={Skamarock, William C and Klemp, Joseph B and Duda, Michael G and Fowler, Laura D and Park, Sang-Hun and Ringler, Todd D},
  journal={Mon.\ Wea.\ Rev.},
  volume={140},
  number={9},
  pages={3090--3105},
  year={2012},
  doi={10.1175/mwr-d-11-00215.1}
}

@article{natale2016compatible,
  title={Compatible finite element spaces for geophysical fluid dynamics},
  author={Natale, Andrea and Shipton, Jemma and Cotter, Colin J},
  journal={Dyn. Stat. Clim. Syst.},
  volume={1},
  number={1},
  pages={1--31},
  year={2016},
  publisher={Oxford University Press UK},
  doi={10.1093/climsys/dzw005}
}

@article{bendall2020compatible,
  title={A compatible finite-element discretisation for the moist compressible {E}uler equations},
  author={Bendall, Thomas M and Gibson, Thomas H and Shipton, Jemma and Cotter, Colin J and Shipway, Ben},
  journal={Quart.\ J.\ Roy.\ Meteor.\ Soc.},
  volume={146},
  number={732},
  pages={3187--3205},
  year={2020},
  publisher={Wiley Online Library},
  doi={10.1002/qj.3841}
}

@article{durran1982effects,
  title={On the effects of moisture on the {B}runt--{V}{\"a}is{\"a}l{\"a} frequency},
  author={Durran, Dale R and Klemp, Joseph B},
  journal={J.\ Atmos.\ Sci.},
  volume={39},
  number={10},
  pages={2152--2158},
  year={1982},
  doi={10.1175/1520-0469(1982)039<2152:OTEOMO>2.0.CO;2}
}

@article{lauritzen2010conservative,
  title={A conservative semi-{L}agrangian multi-tracer transport scheme {(CSLAM)} on the cubed-sphere grid},
  author={Lauritzen, Peter H and Nair, Ramachandran D and Ullrich, Paul A},
  journal={J. Comput. Phys.},
  volume={229},
  number={5},
  pages={1401--1424},
  year={2010},
  publisher={Elsevier},
  doi={10.1016/j.jcp.2009.10.036}
}

@incollection{nair2011emerging,
  title={Emerging numerical methods for atmospheric modeling},
  author={Nair, Ramachandran D and Levy, Michael N and Lauritzen, Peter H},
  booktitle={{N}umerical {T}echniques for {G}lobal {A}tmospheric {M}odels},
  pages={251--311},
  year={2011},
  publisher={Springer},
  doi={10.1007/978-3-642-11640-7_9}
}

@book{gibson2019compatible,
  title={Compatible Finite Element Methods for Geophysical Flows: {A}utomation and Implementation Using {F}iredrake},
  author={Gibson, Thomas H and McRae, Andrew TT and Cotter, Colin J and Mitchell, Lawrence and Ham, David A},
  year={2019},
  publisher={Springer Nature},
  doi={10.1007/978-3-030-23957-2}
}

@manual{FiredrakeUserManual,
  title        = {Firedrake User Manual},
  author       = {David A. Ham and Paul H. J. Kelly and Lawrence Mitchell and Colin J. Cotter and Robert C. Kirby and Koki Sagiyama and Nacime Bouziani and Sophia Vorderwuelbecke and Thomas J. Gregory and Jack Betteridge and Daniel R. Shapero and Reuben W. Nixon-Hill and Connor J. Ward and Patrick E. Farrell and Pablo D. Brubeck and India Marsden and Thomas H. Gibson and Miklós Homolya and Tianjiao Sun and Andrew T. T. McRae and Fabio Luporini and Alastair Gregory and Michael Lange and Simon W. Funke and Florian Rathgeber and Gheorghe-Teodor Bercea and Graham R. Markall},
  organization = {Imperial College London and University of Oxford and Baylor University and University of Washington},
  edition      = {First edition},
  year         = {2023},
  month        = {5},
  doi          = {10.25561/104839},
}

@article{satoh2008nonhydrostatic,
  title={Nonhydrostatic icosahedral atmospheric model ({NICAM}) for global cloud resolving simulations},
  author={Satoh, Masaki and Matsuno, Taro and Tomita, Hirofumi and Miura, Hiroaki and Nasuno, Tomoe and Iga, Shin-Ichi},
  journal={J. Comput. Phys.},
  volume={227},
  number={7},
  pages={3486--3514},
  year={2008},
  publisher={Elsevier},
  doi={10.1016/j.jcp.2007.02.006}
}

@article{melvin2018choice,
  title={Choice of function spaces for thermodynamic variables in mixed finite-element methods},
  author={Melvin, Thomas and Benacchio, Tommaso and Thuburn, John and Cotter, Colin},
  journal={Quart.\ J.\ Roy.\ Meteor.\ Soc.},
  volume={144},
  number={712},
  pages={900--916},
  year={2018},
  publisher={Wiley Online Library},
  doi={10.1002/qj.3268}
}

@article{shu1988efficient,
  title={Efficient implementation of essentially non-oscillatory shock-capturing schemes},
  author={Shu, Chi-Wang and Osher, Stanley},
  journal={J. Comput. Phys.},
  volume={77},
  number={2},
  pages={439--471},
  year={1988},
  publisher={Elsevier},
  doi={10.1016/0021-9991(88)90177-5}
}

@article{wood2014inherently,
  title={An inherently mass-conserving semi-implicit semi-{L}agrangian discretization of the deep-atmosphere global non-hydrostatic equations},
  author={Wood, Nigel and Staniforth, Andrew and White, Andy and Allen, Thomas and Diamantakis, Michail and Gross, Markus and Melvin, Thomas and Smith, Chris and Vosper, Simon and Zerroukat, Mohamed and others},
  journal={Quart.\ J.\ Roy.\ Meteor.\ Soc.},
  volume={140},
  number={682},
  pages={1505--1520},
  year={2014},
  publisher={Wiley Online Library},
  doi={10.1002/qj.2235}
}

@article{kuzmin2010vertex,
  title={A vertex-based hierarchical slope limiter for p-adaptive discontinuous {G}alerkin methods},
  author={Kuzmin, Dmitri},
  journal={J. Comput. Appl. Math.},
  volume={233},
  number={12},
  pages={3077--3085},
  year={2010},
  publisher={Elsevier},
  doi={10.1016/j.cam.2009.05.028}
}

@article{bosler2019conservative,
  title={Conservative multimoment transport along characteristics for discontinuous {G}alerkin methods},
  author={Bosler, Peter A and Bradley, Andrew M and Taylor, Mark A},
  journal={SIAM J. Sci. Comput.},
  volume={41},
  number={4},
  pages={B870--B902},
  year={2019},
  publisher={SIAM},
  doi={10.1137/18m1165943}
}

@article{giraldo1997lagrange,
  title={Lagrange--{G}alerkin methods on spherical geodesic grids},
  author={Giraldo, Francis X},
  journal={J. Comput. Phys.},
  volume={136},
  number={1},
  pages={197--213},
  year={1997},
  publisher={Elsevier},
  doi={10.1006/jcph.1997.5771}
}

@article{nair2005discontinuous,
  title={A discontinuous {G}alerkin transport scheme on the cubed sphere},
  author={Nair, Ramachandran D and Thomas, Stephen J and Loft, Richard D},
  journal={Mon.\ Wea.\ Rev.},
  volume={133},
  number={4},
  pages={814--828},
  year={2005},
  doi={10.1175/mwr2890.1}
}

@incollection{machenhauer2009finite,
  title={Finite-volume methods in meteorology},
  author={Machenhauer, Bennert and Kaas, Eigil and Lauritzen, Peter Hjort},
  booktitle={Handbook of Numerical Analysis},
  volume={14},
  pages={3--120},
  year={2009},
  publisher={Elsevier},
  doi={10.1016/s1570-8659(08)00201-9}
}

@article{mcrae2016automated,
  title={Automated generation and symbolic manipulation of tensor product finite elements},
  author={McRae, Andrew TT and Bercea, G-T and Mitchell, Lawrence and Ham, David A and Cotter, Colin J},
  journal={SIAM Journal on Scientific Computing},
  volume={38},
  number={5},
  pages={S25--S47},
  year={2016},
  publisher={SIAM},
  doi={10.1137/15m1021167}
}

@article{lauritzen2012evaluating,
  title={Evaluating advection/transport schemes using interrelated tracers, scatter plots and numerical mixing diagnostics},
  author={Lauritzen, Peter Hjort and Thuburn, John},
  journal={Quarterly Journal of the Royal Meteorological Society},
  volume={138},
  number={665},
  pages={906--918},
  year={2012},
  publisher={Wiley Online Library},
  doi={10.1002/qj.986}
}

@article{lauritzen2012standard,
  title={A standard test case suite for two-dimensional linear transport on the sphere},
  author={Lauritzen, Peter H and Skamarock, William C and Prather, Michael J and Taylor, Mark A},
  journal={Geosci. Model Dev.},
  volume={5},
  number={3},
  pages={887--901},
  year={2012},
  doi={10.5194/gmd-5-887-2012}
}

@article{zhang2010maximum,
  title   = {On maximum-principle-satisfying high order schemes for scalar conservation laws},
  author  = {Zhang, Xiangxiong and Shu, Chi-Wang},
  journal = {Journal of Computational Physics},
  volume  = {229},
  number  = {9},
  pages   = {3091--3120},
  year    = {2010},
  doi     = {10.1016/j.jcp.2009.12.030}
}

@article{kuzmin2010failsafe,
  title   = {Failsafe flux limiting and constrained data projections for equations of gas dynamics},
  author  = {Kuzmin, Dmitri and M{\"o}ller, Matthias and Shadid, John N. and Shashkov, Mikhail},
  journal = {Journal of Computational Physics},
  volume  = {229},
  number  = {23},
  pages   = {8766--8779},
  year    = {2010},
  doi     = {10.1016/j.jcp.2010.08.009}
}

@article{bradley2019communication,
  title   = {Communication-efficient property preservation in tracer transport},
  author  = {Bradley, Andrew M. and Bosler, Peter A. and Guba, Oksana and Taylor, Mark A. and Barnett, Gregory A.},
  journal = {SIAM Journal on Scientific Computing},
  volume  = {41},
  number  = {3},
  pages   = {C161--C193},
  year    = {2019},
  doi     = {10.1137/18M1165414}
}

@misc{shipton2026gusto,
  title         = {Gusto: a toolkit for compatible finite element dynamical cores},
  author        = {Jemma Shipton and Thomas M. Bendall and Colin J. Cotter and David M. Acreman and Timothy C. Andrews and Jack Betteridge and Alex Brown and Thomas H. Gibson and Nell Hartney and David A. Ham and Joshua Hope-Collins and Jeremy A. Pike and Connor J. Ward and Dan Witt},
  year          = {2026},
  eprint        = {2611.17540},
  archivePrefix = {arXiv},
  primaryClass  = {math.NA},
  doi           = {10.48550/arXiv.2611.17540},
  url           = {https://arxiv.org/abs/2611.17540},
  note          = {Preprint}
}

@article{konor2000choice,
  title={Choice of a vertical grid in incorporating condensation heating into an isentropic vertical coordinate model},
  author={Konor, Celal S and Arakawa, Akio},
  journal={Mon.\ Wea.\ Rev.},
  volume={128},
  number={11},
  pages={3901--3910},
  year={2000},
  publisher={American Meteorological Society},
  doi={10.1175/1520-0493(2001)129<3901:COAVGI>2.0.CO;2}
}

@article{thuburn2005vertical,
  author  = {Thuburn, J. and Woollings, T. J.},
  title   = {Vertical discretizations for compressible {E}uler equation atmospheric models giving optimal representation of normal modes},
  journal = {Journal of Computational Physics},
  volume  = {203},
  number  = {2},
  pages   = {386--404},
  year    = {2005},
  doi     = {10.1016/j.jcp.2004.08.018}
}

@article{lee2025helmholtz,
  author  = {Lee, David and Mart{\'i}n, Alberto F. and Ricardo, Kieran},
  title   = {Helmholtz preconditioning for the compressible {E}uler equations using mixed finite elements with {L}orenz staggering},
  journal = {Quarterly Journal of the Royal Meteorological Society},
  volume  = {151},
  number  = {766},
  pages   = {e4887},
  year    = {2025},
  doi     = {10.1002/qj.4887}
}

@article{guba2014optimization,
  title   = {Optimization-based limiters for the spectral element method},
  author  = {Guba, Oksana and Taylor, Mark and St-Cyr, Amik},
  journal = {J. Comput. Phys.},
  volume  = {267},
  pages   = {176--195},
  year    = {2014},
  doi     = {10.1016/j.jcp.2014.02.029}
}

@article{schar1996synchronous,
  title   = {A synchronous and iterative flux-correction formalism for coupled transport equations},
  author  = {Sch{\"a}r, Christoph and Smolarkiewicz, Piotr K.},
  journal = {J. Comput. Phys.},
  volume  = {128},
  number  = {1},
  pages   = {101--120},
  year    = {1996},
  doi     = {10.1006/jcph.1996.0198}
}

@inproceedings{barth1989design,
  title     = {The design and application of upwind schemes on unstructured meshes},
  author    = {Barth, Timothy J. and Jespersen, Dennis C.},
  booktitle = {27th Aerospace Sciences Meeting},
  note      = {AIAA Paper 89-0366},
  address   = {Reno, NV},
  year      = {1989},
  doi       = {10.2514/6.1989-366}
}

\end{document}